\PassOptionsToPackage{usenames,dvipsnames,x11names}{xcolor}
\documentclass[hidelinks,onefignum,onetabnum,onethmnum,oneeqnum]{siamart220329}

\usepackage{url}
\usepackage{mathtools}
\usepackage{upgreek}
\usepackage{amssymb}
\usepackage{paralist}
\usepackage[shortlabels]{enumitem}
\usepackage[usenames,dvipsnames,x11names]{xcolor}

\usepackage[format=hang,font=small,labelfont={rm,sc},textfont=sl,width=0.95\textwidth]{caption}
\usepackage{subcaption}
\usepackage[bottom]{footmisc}
\usepackage{graphicx}
\newcommand{\brvR}{\brv\Rsup}
\newcommand{\brvRt}{\brv\Rsup{}\tsup}
\newcommand{\bfCr}{\brC}

\newcommand{\bryU}{\bry_{U}}
\newcommand{\bryO}{\bry_{\OO}}
\newcommand{\brfR}{\brf\Rsup}
\newcommand{\brfRt}{\brf\Rsup{}\tsup}
\newcommand{\Tsubr}{_{\text{\tiny T},r}}
\newcommand{\PcalH}{\widehat{\Pcal}}

\newcommand{\brZUU}{\brZ_{\scriptscriptstyle UU}}
\newcommand{\brZUO}{\brZ_{\scriptscriptstyle U\Omega}}
\newcommand{\brZOU}{\brZ_{\scriptscriptstyle \Omega U}}
\newcommand{\brZOO}{\brZ_{\scriptscriptstyle \Omega\Omega}}
\newcommand{\AUU}{A_{\scriptscriptstyle UU}}
\newcommand{\AUO}{A_{\scriptscriptstyle U\Omega}}
\newcommand{\AOO}{A_{\scriptscriptstyle \Omega\Omega}}
\newcommand{\BUU}{B_{\scriptscriptstyle UU}}
\newcommand{\BUO}{B_{\scriptscriptstyle U\Omega}}
\newcommand{\BOO}{B_{\scriptscriptstyle \Omega\Omega}}
\newcommand{\bsfHs}{\bsfH\Tsub}
\newcommand{\bsfHsr}{\bsfH\Tsubr}
\newcommand{\Newcommand}[2]{\providecommand{#1}{}\renewcommand{#1}{#2}}

\makeatletter
\newcommand{\proofstep}[1]{%
  \par% ensure starting on a new paragraph
  \addvspace{\smallskipamount}% some vertical space
  \textit{#1\@addpunct{.}}\enspace\ignorespaces
}
\newlength{\parindentsave}
\newcommand{\Nproofstep}[1]{%
\setlength{\parindent}{0pt}
  \par% ensure starting on a new paragraph
  \addvspace{\smallskipamount}% some vertical space
  \textit{\noindent #1\@addpunct{.}}\enspace\ignorespaces
\setlength{\parindent}{\parindentsave}
}
\makeatother
\newcommand{\bfA} {\boldsymbol{A}}
\newcommand{\bfal}{\boldsymbol{\alpha}}
\newcommand{\bfbe}{\boldsymbol{\beta}}
\newcommand{\bfC} {\boldsymbol{C}}
\newcommand{\bfd} {\boldsymbol{d}}
\newcommand{\bfD} {\boldsymbol{D}}
\newcommand{\bfe} {\boldsymbol{e}}
\newcommand{\bff} {\boldsymbol{f}}
\newcommand{\bffR} {\RBmath{f}}
\newcommand{\bfg} {\boldsymbol{g}}
\newcommand{\bfI} {\boldsymbol{I}}
\newcommand{\bfmu}{\boldsymbol{\mu}}
\newcommand{\bfn} {\boldsymbol{n}}
\newcommand{\bfna}{\boldsymbol{\nabla}}
\newcommand{\bfOm}{\boldsymbol{\Omega}}
\newcommand{\bfPi}{\boldsymbol{\Pi}}
\newcommand{\bfr} {\boldsymbol{r}}
\newcommand{\bfu} {\boldsymbol{u}}
\newcommand{\bfU} {\boldsymbol{U}}
\newcommand{\bfuD} {\bfu\Dsup}
\newcommand{\bfups}{\boldsymbol{\upsilon}}
\newcommand{\bfuR} {\RBmath{u}}
\newcommand{\bfuS} {\bfu\Ssup}
\newcommand{\bfv} {\boldsymbol{v}}
\newcommand{\bfV} {\boldsymbol{V}}
\newcommand{\bfw} {\boldsymbol{w}}
\newcommand{\bfW} {\boldsymbol{W}}
\newcommand{\bfwR} {\RBmath{w}}
\newcommand{\bfx} {\boldsymbol{x}}
\newcommand{\bfy} {\boldsymbol{y}}
\newcommand{\bfz} {\boldsymbol{z}}
\newcommand{\bfZ} {\boldsymbol{Z}}
\newcommand{\bfze}{\mathbf{0}}
\newcommand{\brA} {\mathbf{A}}
\newcommand{\brC} {\mathbf{C}}
\newcommand{\brf} {\mathbf{f}}
\newcommand{\brI} {\mathbf{I}}
\newcommand{\brO} {\mathbf{O}}
\newcommand{\brr} {\mathbf{r}}
\newcommand{\brR} {\mathbf{R}}
\newcommand{\brs} {\mathbf{s}}
\newcommand{\brt} {\mathbf{t}}
\newcommand{\brv} {\mathbf{v}}
\newcommand{\bry} {\mathbf{y}}
\newcommand{\brz} {\mathbf{z}}
\newcommand{\brZ} {\mathbf{Z}}
\newcommand{\bsfH} {\boldsymbol{\mathsf{H}}}
\newcommand{\bsfN} {\boldsymbol{\mathsf{N}}}
\newcommand{\bsfR} {\boldsymbol{\mathsf{R}}}
\newcommand{\del}[1][]{\partial_{#1}}
\newcommand{\dip} {\! :\!}
\newcommand{\Div}{\mbox{div}\,}
\newcommand{\dO}{\partial\OO}
\newcommand{\dotp}{\raisebox{1pt}{\hspace*{1pt}\scalebox{0.45}{$\bullet$}}\hspace*{1pt}}
\newcommand{\dS}{\,\text{d}S}
\newcommand{\Dsub}{_{\text{\tiny D}}}
\newcommand{\Dsup}{^{\text{\tiny D}}}
\newcommand{\dV}{\;\text{d}V}
\newcommand{\eps}{\varepsilon}
\newcommand{\G}{\Gamma}
\newcommand{\iG} {\int_{\G}}
\newcommand{\inv}[1]{\dfrac{1}{#1}}
\newcommand{\iO}{\int_{\OO}}
\newcommand{\labs}{\big|\hspace*{0.1em}}
\newcommand{\lbra}{\big\langle\hspace*{0.1em}}
\newcommand{\lcb}{\big\{\hspace*{0.1em}}
\newcommand{\lpar}{\big(\hspace*{0.1em}}
\newcommand{\Lpar}{\Big(\,}
\newcommand{\lsqb}{\big[\hspace*{0.1em}}
\newcommand{\Lsqb}{\Big[\hspace*{0.1em}}
\newcommand{\OO}{\Omega}
\newcommand{\Pcal}{\mathcal{P}}
\newcommand{\rabs}{\hspace*{0.1em}\big|}
\newcommand{\Rbb} {\mathbb{R}}
\newcommand{\rbra}{\hspace*{0.1em}\big\rangle}
\newcommand{\rcb}{\hspace*{0.1em}\big\}}
\newcommand{\rmA} {\mathrm{A}}
\newcommand{\rmC} {\mathrm{C}}
\newcommand{\rmm} {\mathrm{m}}
\newcommand{\rpar}{\hspace*{0.1em}\big)}
\newcommand{\Rpar}{\,\Big)}
\newcommand{\rsqb}{\hspace*{0.1em}\big]}
\newcommand{\Rsqb}{\hspace*{0.1em}\Big]}
\newcommand{\Rsup}{^{\text{\tiny R}}}
\newcommand{\sfR} {\mathsf{R}}
\newcommand{\shcap}{\hspace*{-0.1em}\cap\hspace*{-0.1em}}
\newcommand{\shcup}{\hspace*{-0.1em}\cup\hspace*{-0.1em}}
\newcommand{\shdeq}{\hspace*{-0.1em}:=\hspace*{-0.1em}}
\newcommand{\sheq}{\hspace*{-0.1em}=\hspace*{-0.1em}}
\newcommand{\shin}{\hspace*{-0.1em}\in\hspace*{-0.1em}}
\newcommand{\shl}{\hspace*{-0.1em}<\hspace*{-0.1em}}
\newcommand{\shleq}{\hspace*{-0.1em}\leq\hspace*{-0.1em}}
\newcommand{\shm}{\hspace*{-0.1em}-\hspace*{-0.1em}}
\newcommand{\shneq}{\hspace*{-0.1em}\not=\hspace*{-0.1em}}
\newcommand{\shp}{\hspace*{-0.1em}+\hspace*{-0.1em}}
\newcommand{\shsetm}{\hspace*{-0.1em}\setminus\hspace*{-0.1em}}
\newcommand{\shsubs}{\hspace*{-0.1em}\subset\hspace*{-0.1em}}
\newcommand{\shtimes}{\hspace*{-0.1em}\times\hspace*{-0.1em}}
\newcommand{\sip} {\! \cdot\!}
\newcommand{\Span}{\text{span}}
\newcommand{\ssub}{_{\text{s}}}
\newcommand{\Ssub}{_{\text{\tiny S}}}
\newcommand{\Ssup}{^{\text{\tiny S}}}
\newcommand{\suite}[1][0ex]{\notag \\[#1] & \mbox{}\hspace{15pt}}
\newcommand{\sym}{_{\text{sym}}}
\newcommand{\Tsub}{_{\text{\tiny T}}}
\newcommand{\tsup}{^{\text{t}}}
\newcommand{\Tsup}{^{\text{\tiny T}}}
\Newcommand{\bfups}{\boldsymbol{\upupsilon}}
\Newcommand{\uS}{u\Ssup}
\Newcommand{\uSp}{(u\Ssup)'}
\Newcommand{\bffR} {\bff\Rsup}
\Newcommand{\bffHR} {\bffH\Rsup}
\Newcommand{\bfuHR} {\bfuH\Rsup}
\Newcommand{\bfuD} {\bfups\Dsup}
\Newcommand{\bfuS} {\bfups\Ssup}
\Newcommand{\bfwS} {\bfw\Ssup}
\Newcommand{\bfvS} {\bfv\Ssup}
\Newcommand{\bfzS} {\bfz\Ssup}
\Newcommand{\bfuR} {\bfups\Rsup}
\Newcommand{\bfUR} {\bfU\Rsup}
\Newcommand{\bfCR} {\bfC\Rsub}
\Newcommand{\bfDR} {\bfD\Rsub}
\Newcommand{\bfwR} {\bfw\Rsup}
\newcommand{\bfupal}{\boldsymbol{\upalpha}}

\usepackage{autonum}
\newsiamremark{remark}{Remark}

\headers{Optimal slip profiles for arbitrary 3D swimmers}{M. Bonnet, K. Das, S. Veerapaneni, H. Zhu}

\author{
Marc Bonnet\thanks{POEMS (CNRS, INRIA, ENSTA), ENSTA Paris, 91120 Palaiseau, France
  (\email{mbonnet@ensta.fr})}
\and Kausik Das\thanks{Department of Mathematics, University of Michigan, Ann Arbor, MI, USA (\email{kausik@umich.edu, shravan@umich.edu})}
\and Shravan Veerapaneni\footnotemark[2]
\and Hai Zhu\thanks{Flatiron Institute, Simons Foundation, New York, NY, USA(\email{hzhu@flatironinstitute.org})}
}

\title{Slip optimization on arbitrary 3D microswimmers: a reduced-dimension and boundary-integral framework
}

\begin{document}
\maketitle

\begin{abstract}
  This article presents a computational framework for determining the optimal slip velocity of a microswimmer with arbitrary three-dimensional geometry suspended in a viscous fluid. The objective is to minimize the hydrodynamic power dissipation required to maintain unit speed along the net swimming direction. By exploiting the linearity of the Stokes equations and the Lorentz reciprocal theorem, we derive an explicit linear operator that maps the tangential surface slip velocity to the resulting rigid-body translational and rotational velocities, effectively decoupling the hydrodynamic boundary value problem from the optimization loop. The a priori infinite-dimensional search space for the slip optimization is reduced to the finite dimension $r$ of rigid-body motions by finding an appropriate subspace of the operator's domain. This reduces the PDE-constrained optimization to a low-dimensional programming problem that can be solved at negligible computational cost once the system matrices are assembled. The optimization algorithm requires 2$r$ auxiliary flow problems that are solved numerically using a high-order boundary integral method. We validate the accuracy of the proposed method and present optimal slip profiles and swimming trajectories for a variety of microswimmer shapes. We investigate the effect of some common geometrical symmetries of the swimmer shape on the resulting optimal motion, and in particular present a modified version of the slip optimization algorithm for axisymmetric shapes, where tangential rigid-body velocities may occur.\medskip
\end{abstract}


\section{Introduction}
The hydrodynamics of microswimmers has attracted sus\-tained research interest, motivated both by the biological imperative to understand how microorganisms such as bacteria, ciliates and algae navigate their environments \cite{LaugaPowers2009, ishikawa2024fluid} and by the engineering challenge of designing artificial microswimmers for biomedical applications, such as targeted drug delivery \cite{dreyfus2005, tsang2020roads, LiEtAl2021}. At these scales, inertial effects are negligible, and fluid motion is governed by the Stokes equations. To achieve net displacement, swimmers must employ non-reciprocal surface deformations or slip velocities, a constraint famously encapsulated by the {\em scallop theorem} \cite{Purcell1977}. The squirmer model, originally introduced by Lighthill \cite{Lighthill1952} and adapted for ciliates by Blake \cite{blake:71}, has become the standard theoretical tool for these problems. It replaces the complex individual beating of cilia with a prescribed tangential slip velocity on the boundary of a rigid body \cite{pedley:16}.

A central question in the design of synthetic swimmers is the optimization of this slip velocity to maximize swimming efficiency. Efficiency is typically quantified by comparing the power dissipated by the active slip to the power required to tow the passive rigid body at the same speed. For spherical geometries, analytical solutions for optimal slip profiles have been established using spherical harmonic expansions \cite{MichelinLauga2010}. More recently, computational approaches have extended these results to axisymmetric shapes, determining optimal gaits for fixed shapes \cite{B-2020-04} and even performing joint optimization of shape and slip \cite{B-2023-10}. However, these studies have largely been confined to bodies of revolution swimming in straight lines along prescribed directions.

The restriction to axisymmetry represents a significant gap in the current literature. Many biological organisms and synthetic particles (e.g., chiral swimmers or general Janus particles) possess arbitrary three-dimensional shapes \cite{samatas:23}. While a minimum dissipation theorem \cite{nasouri:21} establishes that the optimal slip for any shape is identical to the flow around the corresponding perfect-slip body, an efficient computational realization of this result for general three-dimensional geometries has remained lacking. Unlike their axisymmetric counterparts, general 3D swimmers exhibit a complex coupling between translational and rotational velocities, generically resulting in helical trajectories  \cite{Crenshaw1993, DeSimone2020}. Optimizing the slip profile for such shapes is mathematically and computationally demanding. The search space for the surface slip is infinite-dimensional, and determining the resulting rigid-body motion---required to evaluate the objective function---necessitates solving the forward Stokes problem. In a naive optimization loop, this would incur a prohibitive computational cost, as the flow field must be resolved for every trial slip profile to satisfy the force-free and torque-free conditions.

In this work, we present a fast and rigorous computational framework for determining the optimal slip profiles of microswimmers with arbitrary 3D shapes. Our approach relies on the linearity of the Stokes equations and the Lorentz reciprocal theorem to decouple the flow solution from the optimization variables. We derive a linear operator that maps the tangential slip velocity directly to the resulting swimming velocity via auxiliary ``extractor'' solutions. This allows us to strictly reduce the infinite-dimensional search space to a finite-dimensional subspace, which depends only on the swimmer shape and whose dimension is that of the space of rigid-body velocities on the swimmer surface, capable of generating any given rigid-body motion at least power loss; the value of the latter is shown to be that established in~\cite{nasouri:21}. Consequently, efficiency optimization is formulated as a low-dimensional programming problem that can be solved rapidly without repeated flow solves. We utilize a high-order boundary integral method (BIM) to discretize the problem, which handles the unbounded fluid domain exactly and requires meshing only the swimmer surface. While the primary focus of this article is on formulation and algorithmic development, our companion paper~\cite{B-2025-05} analyzes the physical implications, studies the helical trajectories, and derives some analytical results for spheroidal microswimmers.

The paper is organized as follows. In Section~\ref{sec:forward:opt}, we define the forward problem for a swimmer of arbitrary shape and state the power loss minimization problem. Section~\ref{sec:rdim} details the construction of the finite-dimensional slip velocity search space using auxiliary flow solutions, which allows for the exact reduction of the problem dimension. We also introduce and validate the boundary integral formulation used to obtain the auxiliary flow solutions. In Section~\ref{sec:solv:opt}, we develop the solution method for the resulting low-dimensional optimization problem to determine the optimal gait, taking advantage of a further simplification permitted by the quadratic character of the partial minimization with fixed swimming direction. Section~\ref{CA:sym} addresses several common cases of geometrical symmetry in the swimmer shape, and the slip optimization for axisymmetric shapes is fully treated in Section~\ref{sec:axi}. Section~\ref{sec:conclu} concludes the paper. Numerical validations of the main steps are provided in Sections~\ref{sec:rdim}, \ref{sec:solv:opt} and~\ref{sec:axi}.

\section{Forward and optimization problems}
\label{sec:forward:opt}

\subsection{Forward flow problem}

Let the microswimmer occupy the $d$-dimensional bounded domain $\OO\Ssub\shsubs\Rbb^d$ with (closed, smooth) boundary $\dO\ssub=\G$, and let $\Omega=\Rbb^d\shsetm\overline{\OO\Ssub}$ denote the unbounded fluid region surrounding it. We will let $\bfn$ denote the unit normal on $\G$ directed away from the fluid. In the low Reynolds number limit, fluid flows are assumed to obey the incompressible Stokes equations
\begin{equation}
  -\bfna p + \mu\Delta\bfu = \bfze, \quad \Div\bfu = 0 \qquad \text{in $\OO$} \label{Stokes}
\end{equation}
(where $\bfu$ is the velocity field, $p$ the pressure and $\mu$ the dynamic viscosity) and to arise from prescribing the velocity on $\G$ and a decay condition at infinity:
\begin{equation}
  \bfu = \bfuD \quad \text{on $\G$}, \qquad \lim_{|\bfx|\to\infty} \bfu(\bfx)=\bfze. \label{BC}
\end{equation}
The datum $\bfuD$ must satisfy the compatibility condition $\lbra\bfuD,\bfn\rbra_{\G}=0$, where $\lbra\cdot,\cdot\rbra_{\G}$ denotes the $L^2(\G)$ scalar product (with pointwise dot or matrix-vector products implied when used for vector- or matrix-valued functions). Problem~\eqref{Stokes}-\eqref{BC} is well-posed for $\bfuD\in\bsfH:=\lcb \bfw\shin H^{1/2}(\G;\Rbb^d),\, \lbra\bfw,\bfn\rbra_{\G}\sheq 0 \rcb$, see e.g.~\cite[Sec.~I.2]{temam}. In what follows, the traction field $\bff:= -p\bfn\shp2\mu\bfD[\bfu]\sip\bfn$ on $\G$, where $\bfD[\bfu]:=\tfrac{1}{2}(\bfna\bfv\shp\bfna\tsup\bfv)$ denotes the strain rate tensor, will play an important role; its dependence on the prescribed velocity $\bfuD$ will be emphasized, when needed, by the notation $\bff[\bfuD]$. In general, $\bff[\bfuD]$ has non-zero net forces and torques.\enlargethispage*{3ex}

For any Dirichlet datum $\bfuD$, problem~\eqref{Stokes}-\eqref{BC} obeys the following well-known sign reciprocity properties, see e.g.~\cite[Sec.~3]{masoud:stone:19}, which will play an important role; they are briefly proved for completeness in Sec.~\ref{proof:recipr}.
\begin{lemma}[Lorentz reciprocity identity]\label{recipr}
Any solution of problem~\eqref{Stokes}-\eqref{BC}, with datum $\bfuD$, verifies
\begin{equation}
  \lbra \bfuD,\bff[\bfuD] \rbra_{\G} = \iO 2\mu\bfD[\bfu]\dip\bfD[\bfu] \dV.
\end{equation}
Any pair of solutions of problem~\eqref{Stokes}-\eqref{BC}, with respective data $\bfuD_1$ and $\bfuD_2$, verifies
\begin{equation}
  \lbra \bfuD_1,\bff[\bfuD_2] \rbra_{\G} = \lbra \bfuD_2,\bff[\bfuD_1] \rbra_{\G};
\end{equation}
Consequently, the bilinear form $(\bfuD_1,\bfuD_2)\mapsto \lbra \bfuD_1,\bff[\bfuD_2] \rbra_{\G}$ is symmetric and positive.
\end{lemma}

For swimmers, the velocity datum $\bfuD$ has the form
\begin{equation}
  \bfuD = \bfuS + \bfuR, \qquad \bfuS\shin\bsfHs,\ \bfuR=\bfU + \bfOm\shtimes\bfx\in\bsfR,
\label{BC:swimmer}
\end{equation}
where $\bfuS$ is a given (tangential) slip velocity, $\bsfHs = \lcb \bfw\in\bsfH,\, \bfw\tsup\bfn=0 \rcb$ is the subspace of $\bsfH$ comprised of tangential velocity fields, $\bsfR$ is the finite-dimensional subspace of $\bsfH$ containing all traces on $\G$ of rigid-body velocities (whose dimension will be denoted as $r:=\text{dim}(\bsfR)$ throughout), and the swimming velocity $\bfuR$ (i.e the value of $\bfU,\bfOm$) is determined for given $\bfuS$ by additionally requiring that the traction $\bff[\bfuD]=\bff[\bfuS]\shp\bff[\bfuR]$ have zero net force and net torque, i.e.
\begin{equation}
  \lbra \bff[\bfuD], \bfwR \rbra_{\G} = 0 \qquad \text{for all } \bfwR\shin\bsfR. \label{no-net-force}
\end{equation}
In this work, we focus on three-dimensional swimmers and flows, i.e. $d\sheq3$ (for which $r\sheq6$), while $r\sheq3$ or $r\sheq1$ for, respectively, two-dimensional or axisymmetric flows.

For any given slip profile $\bfuS$ that is time-independent in the body frame, the induced rigid-body velocity $\bfU[\bfuS]\shp\bfOm[\bfuS]\shtimes\bfx$ is shown in~\cite{B-2025-05} to result in the swimmer centroid following a helical trajectory whose translation velocity $\bfW=\bfW[\bfuS]$ along the helix axis is given by
\begin{equation}
  \text{(a) \ }\bfW
 = \frac{\bfU\tsup\bfOm}{\labs\bfOm\rabs^2} \, \bfOm \quad (\text{if }\bfOm\shneq\bfze), \qquad
  \text{(b) \ }\bfW = \bfU \quad (\text{if }\bfU\shtimes\bfOm\sheq\bfze). \label{W:def}
\end{equation}
Case~(\ref{W:def}a) covers all possible motions such that $\bfOm\not=\bfze$, which are either circular motions (if $\bfU\tsup\bfOm\sheq 0$, i.e. $\bfW\sheq\bfze$), spinning straight motions (if $\bfU\shneq\bfze$, $\bfU\shtimes\bfOm\sheq\bfze$) or helical motions (otherwise). The net velocity $\bfW$ as given by~(\ref{W:def}a) does not depend on the rotation velocity magnitude $|\bfOm|$ (whereas the incurred power loss \emph{a priori} does), and is not defined for $\bfOm\sheq\bfze$. The special case~(\ref{W:def}b) covers all straight (i.e. non-helical) motions, which may be either pure translations ($\bfOm\sheq\bfze$) or straight spinning motions ($\bfU\shtimes\bfOm\sheq\bfze,\,\bfOm\shneq\bfze$); it then coincides with case~(\ref{W:def}a) if $\bfOm\not=\bfze$ while allowing for $\bfOm\sheq\bfze$.

From~\eqref{W:def}, all rigid-body motions producing a nonzero net translation velocity $\bfW$ have the form
\begin{equation}
  \bfU = \bfW\shp\bfV, \quad \bfOm = s\bfW, \qquad
  \text{with \ } s\shin\Rbb,\ \bfV\shin\Rbb^d,\ \bfW\tsup\bfV=0, \label{rigid:param}
\end{equation}
and may therefore be parametrized using $(\bfW,s,\bfV)$ instead of $\bfU,\bfOm$. This parametrization is consistent with~(\ref{W:def}a) if $s\shneq0$, and with~(\ref{W:def}b) if $s\shin\Rbb$ and $\bfV\sheq\bfze$; on the other hand, \eqref{rigid:param} with $s\sheq0,\,\bfV\shneq\bfze$ is inconsistent with~(\ref{W:def}b). When $s\shneq0$ and $\bfV\shneq\bfze$, the rigid-body motion follows a helix of radius $|\bfV|/|s\bfW|$ and pitch $2\pi/|s|$, the circular component of the motion having a period $2\pi/|s\bfW|$.

\subsection{Evaluation of swimming velocities}

We now set the framework and notations pertaining to the arising rigid-body motions. Any rigid-body velocity field $\bfuR\shin\bsfR$ may be set in the form
\begin{equation}
  \bfuR(\bfx) = \brvR(\bfx)\bfal, \quad \bfx\shin\G \label{uR:exp}
\end{equation}
where the $d\shtimes r$ matrix-valued function $\brvR=\lsqb \bfuR_1,\, \ldots,\, \bfuR_r \rsqb$ collects $r$ basis vector functions arranged in columns (this notation will similarly be used for other finite sets of vector fields on $\G$) and $\bfal\shin\Rbb^r$ is a (column) vector of scalar coefficients. For the general 3D case, we use $\bfuR_i=\bfe_i$ and $\bfuR_{i+3}=\bfe_i\shtimes\bfx=\eps_{iqp}x_p\bfe_q$ for $i=1,2,3$ as basis functions, so that $\bfal = [\bfU; \bfOm]$ with $\bfU$ and $\bfOm$ collecting the Cartesian components of the rigid-body translation and rotation velocities. Then, for each basis function $\bfuR_{\ell}$, let $\bffR_{\ell}:=\bff[\bfuR_{\ell}]$ be the traction field of the no-slip flow solving problem~\eqref{Stokes}-\eqref{BC} with $\bfuD=\bfuR_{\ell}|_{\G}$. The traction $\bff[\bfuS\shp\bfuR]$, with $\bfuR$ of the form~\eqref{uR:exp}, is thus given by
\begin{equation}
  \bff[\bfuD](\bfx) = \bff[\bfuS](\bfx) + \brfR(\bfx)\bfal, \quad \bfx\shin\G \label{f:xi}
\end{equation}
having set $\brfR(\bfx) := [\bffR_1(\bfx),\, \ldots,\, \bffR_r(\bfx)]\in\Rbb^{d\times r}$. The no-net-force conditions~\eqref{no-net-force} then reduce to finding $\bfal=\bfupal[\bfuS]$ satisfying
\begin{equation}
  \lbra \brvRt,\bff[\bfuS]+\brfR\bfal \rbra_{\G}
  = \lbra \brfRt,\bfuS \rbra_{\G} + \bfCr\bfal = \bfze \label{CR:eq}
\end{equation}
where the matrix $\bfCr\in\Rbb^{r\times r}\sym$ is given by
\begin{equation}
  \bfCr := \lbra \brfRt,\brvR \rbra_{\G} = \lbra \brvRt,\brfR \rbra_{\G}, \quad\text{i.e.}\quad
  \bfCr_{k\ell} = \lbra \bffR_k,\bfuR_{\ell} \rbra_{\G} = \lbra \bfuR_k,\bffR_{\ell} \rbra_{\G} = \bfCr_{\ell k}. \label{CR:expr}
\end{equation}
Lemma~\ref{recipr}, from which in particular the above symmetry of $\bfCr$ follows, is used in both~\eqref{CR:eq} and~\eqref{CR:expr}. Solving~\eqref{CR:eq} for $\bfal$, each entry $\alpha_i[\bfuS]$ of $\bfupal[\bfuS]$ is obtained in the form of a power integral:
\begin{equation}
  \alpha_i[\bfuS] = \lbra \bff^{\alpha}_i,\bfuS \rbra_{\G}, \qquad\text{with}\quad
  \bff^{\alpha}_i := -(\brfR\bfCr^{-1})_i, \label{R:expr}
\end{equation}
with the traction fields $\bff^{\alpha}_i$ acting as the ``extractors'' of the components of $\bfal$. Equation~\eqref{R:expr} defines a linear operator $\sfR:\bsfHs\to\bsfR$ such that $\bfuR=\sfR\bfuS$ is the swimming velocity associated through conditions~\eqref{no-net-force} to a given slip velocity $\bfuS\shin\bsfHs$.
\begin{remark}[no-slip resistance tensor]\label{rem:RNS}
Since $\brC\bfal$ is the vector of net hydrodynamic forces and torques induced by applying the rigid-body velocity $\brvR\bfal$ on $\G$, $\brC$ is in fact the no-slip resistance tensor of the swimmer (denoted by $\brR_{\text{NS}}$ in~\cite{nasouri:21}).
\end{remark}

\subsection{Power loss minimization}

We consider slip velocity optimization problems aimed at achieving low-Reynolds locomotion with minimum energy dissipation, for swimmers with given shape $\G$. The power loss for the swimmer motion with given slip velocity $\bfuS$ is
\begin{equation}
  P(\bfuS) := \lbra \bff[\bfuS\shp\sfR\bfuS],\bfuS\shp\sfR\bfuS \rbra_{\G}
  = \lbra \bff[\bfuS\shp\sfR\bfuS],\bfuS \rbra_{\G} \label{Jw:def}
\end{equation}
(the second equality stemming from the no-net-force condition~\eqref{no-net-force}); it is clearly quadratic in $\bfuS$ and, by Lemma~\ref{recipr}, positive.

Our general goal is to find a slip velocity $\bfuS$ that minimizes the power loss $P$ while maintaining the magnitude $|\bfW|$ of the net translation velocity $\bfW[\bfuS]$. Indeed, since $P$ is quadratic (while $\bfW[\bfuS]$ is homogeneous with degree 1) in $\bfuS$, a normalization constraint on $\bfuS$ is necessary. Here we will enforce $|\bfW|=1$ without loss of generality, and our generic optimization problem is:
\begin{equation}
  \min_{\bfuS\in\bsfHs} P(\bfuS) \qquad\text{subject to} \qquad|\bfW[\bfuS]|^2=1. \label{JL:opt:generic}
\end{equation}

Solving problem~\eqref{JL:opt:generic} for $\bfuS$ in turn yields the parameters $\bfU,\bfOm$ of the resulting helical motion and, through~\eqref{W:def}, the net motion direction $\bfW$.
The norm constraint in~\eqref{JL:opt:generic} is quadratic in $\bfuS$ only when $\bfU[\bfuS]$ and $\bfOm[\bfuS]$ are collinear, see~\eqref{W:def}. Problem~\eqref{JL:opt:generic} therefore cannot in general be reduced to a linear eigenvalue problem for $\bfuS$. By contrast, using a comparison power (quadratic in $\bfuS$) for normalization purposes, as done for efficiency optimization, leads to a linear eigenvalue problem; this was exploited in e.g.~\cite{B-2023-10} dealing with axisymmetric swimmers and axial rotationless net motions. 

\begin{remark}
The power loss arising from dragging a rigid body of same shape as the swimmer at the swim velocity $\sfR\bfuS$, given by
\begin{equation}
  J\Dsub(\bfuS) := \lbra \bff[\sfR\bfuS],\sfR\bfuS \rbra_{\G}
 = \bfal\tsup[\bfuS]\,\bfCr\,\bfupal[\bfuS] \label{JD:def}
\end{equation}
(with the second equality resulting from the explicit expression~\eqref{R:expr} of the operator $\sfR$), is quadratic and (by Lemma~\ref{recipr}) positive. For a fixed net motion direction $\bfW$, problem~\eqref{JL:opt:generic} is equivalent to maximizing the efficiency $J\Dsub(\bfuS) /P(\bfuS)$, as done e.g. in~\cite{B-2020-04,B-2023-10}. However, efficiency and power loss minimization become differing goals if the net motion direction $\bfW$ is unknown, and we focus in this work on the latter goal.
\end{remark}

\section{Optimal reduced-dimension representation of slip velocities}
\label{sec:rdim}

The optimization problems~\eqref{JL:opt:generic} \emph{a priori} involve an infinite-dimensional search space. However, the linear mapping $\bsfHs\to\Rbb^r:\bfuS\mapsto\bfupal[\bfuS]$ has a $r$-dimensional range, and thus a null space $\bsfN\shsubs\bsfHs$ of codimension $r$.
Additive decompositions $\bsfHs=\bsfH_r+\bsfN$ of $\bsfHs$ with $\text{dim}(\bsfH_r)=r$, such that
\begin{equation}
  \bfuS = \bfuS_1+\bfuS_2 \quad (\bfuS_1,\bfuS_2) \in \bsfH_r\shtimes\bsfN, \qquad \implies \qquad
  \bfupal[\bfuS]=\bfupal[\bfuS_1]
\end{equation}
for any $\bfuS\shin\bsfHs$, thus exist. Any such subspace $\bsfH_r$ clearly allows to generate all possible swimming rigid-body motions. Moreover, this opens the possibility of defining $\bsfH_r$ so that the low-dimensional optimization problem
\begin{equation}
  \min_{\bfuS\in\bsfH_r} P(\bfuS) \qquad\text{subject to} \qquad|\bfW[\bfuS]|=1. \label{JL:opt:lowdim}
\end{equation}
is equivalent to the original problem~\eqref{JL:opt:generic}, a task for which we propose in what follows a computationally constructive procedure.

\subsection{Construction of $\bsfH_r$}

To make problem~\eqref{JL:opt:lowdim} equivalent to~\eqref{JL:opt:generic}, $\bsfH_r$ should be such that $P(\bfuS_1\shp\bfuS_2)\geq P(\bfuS_1)$ for any $(\bfuS_1,\bfuS_2)\in\bsfH_r\shtimes\bsfN$.
Since the quadratic form $P(\bfuS_1\shp\bfuS_2)$, given by~\eqref{Jw:def}, expands as
\begin{equation}
  P(\bfuS_1\shp\bfuS_2)
 = P(\bfuS_1) + P(\bfuS_2) + 2\lbra \bff[\bfuS_1\shp\sfR\bfuS_1],\bfuS_2 \rbra_{\G} \label{Jw:expand}
\end{equation}
(using Lemma~\ref{recipr} and the assumption $\sfR\bfuS_2=\bfze$), $P(\bfuS_1\shp\bfuS_2)\geq P(\bfuS_1)$ holds for any $(\bfuS_1,\bfuS_2)\in\bsfH_r\shtimes\bsfN$ provided $\lbra \bff[\bfuS_1\shp\sfR\bfuS_1],\bfuS_2 \rbra_{\G}=0$. We now proceed to show that this requirement is met by setting $\bsfH_r=\Span(\bfz_1,\ldots,\bfz_r)$, with the $r$ linearly-independent slip velocities $\bfz_i$ given by
\begin{equation}
  \bfz_i := (\bfv_i\shm\bfv_i\Rsup)|_{\G} \in \bsfHs \label{zi:def:JL}
\end{equation}
in terms of velocity fields $\bfv_i$ and rigid-body velocities $\bfv_i\Rsup$ solving the flow problems
\begin{equation}
\left.\begin{aligned}
  -\bfna p_i + \mu\Delta\bfv_i = \bfze, \qquad \Div\bfv_i = 0 & \qquad \text{in $\OO$} \\
  \bfPi\bff[\bfv_i] = \bfPi\bff^{\alpha}_i, \qquad \lpar \bfv_i\shm\bfv_i\Rsup\rpar\tsup\bfn=0 & \qquad \text{on $\G$} \\
  \lbra \bff[\bfv_i],\bfwR \rbra_{\G} = 0 & \qquad \text{for all }\bfwR\shin\bsfR \\
  \lim_{|\bfx|\to\infty} \bfv_i(\bfx)=\bfze &
\end{aligned}\ \right\}\qquad (1\shleq i\shleq r), \label{forward:PDE:fD:JL}
\end{equation}
where $\bfPi(\bfx):=\bfI\shm\bfn(\bfx)\!\otimes\!\bfn(\bfx)$ is the orthogonal projector on the tangent plane at $\bfx\shin\G$ and $\bff^{\alpha}_i$ is the traction field introduced as slip extractor in~\eqref{R:expr}.

The flow problems~\eqref{forward:PDE:fD:JL} are designed so that, for each $i\in\{1,\ldots,r\}$ (i) $\bfz_i\shdeq(\bfv_i\shm\bfv\Rsup_i)|_{\G}$ is a tangential (slip) velocity and (ii) the field $\bfv_i$ develops the tangential tractions allowing the rigid-body extraction of $\alpha_i[\bfuS]$ via~\eqref{R:expr}.
The flow fields $(\bfv_i,p_i,\bfv_i\Rsup)$ and associated slip velocities $\bfz_i$ have the following properties (see proof in Sec.~\ref{lemma:z:proof}):
\begin{lemma}\label{lemma:z}
Let $\bsfH_r\shdeq\Span(\bfz_1,\ldots,\bfz_r)$ with the slip velocities $\bfz_i$ defined by~\eqref{zi:def:JL} in terms of the respective solutions of problems~\eqref{forward:PDE:fD:JL}. Let $\brz(\bfx)=[\bfz_1(\bfx),\ldots,\bfz_r(\bfx)]$.
\begin{compactenum}[(i)]
\item The solutions of problems~\eqref{forward:PDE:fD:JL} solve problems~\eqref{Stokes}-\eqref{BC} with $\bfuD=\bfz_i\shp\bfv_i\Rsup$ and verify the no-net-force conditions~\eqref{no-net-force}; in particular, we have $\bfv\Rsup_i=\sfR\bfz_i$.
\item For any $\bfw\shin\bsfHs$, each $\bfz_i$ verifies
$\lbra \bff[\bfz_i\shp\sfR\bfz_i],\bfw \rbra_{\G} = \lbra \bff[\bfv_i],\bfw \rbra_{\G} = \alpha_i[\bfw].$
\item Let $\bfuS\shin\bsfHs$ and write $\bfuS=\bfuS_1\shp\bfuS_2$ with $(\bfuS_1,\bfuS_2)\in\bsfH_r\shtimes\bsfN$. We have
\begin{equation}
  P(\bfuS_1\shp\bfuS_2) = P(\bfuS_1) + P(\bfuS_2) \geq P(\bfuS_1).  \label{Jw:reduc}
\end{equation}
\end{compactenum}
\end{lemma}

We next use the $\bfv_i$ to define the matrix $\brA=[\rmA_{ij}]\shin\Rbb^{r\times r}$ given by
\begin{equation}
  \rmA_{ij}
 = \lbra \bff[\bfv_i],\bfv_j \rbra_{\G}
 = \lbra \bff[\bfz_i\shp\sfR\bfz_i],\bfz_j \rbra_{\G} \quad 1\shleq i,j\shleq r \label{eig:rank:r:JL}
\end{equation}
(with the second equality stemming from~\eqref{zi:def:JL} and the no-net-force condition in problem~\eqref{forward:PDE:fD:JL}). The matrix $\brA$ is symmetric and positive definite by Lemma~\ref{recipr}; moreover, property (ii) of Lemma~\ref{lemma:z} implies
\begin{equation}
  \rmA_{ij} = \alpha_i[\bfz_j]. \label{A:zi}
\end{equation}
We then introduce the set $\bry:=[\bfy_1,\, \ldots,\, \bfy_r]$ of tangential slip velocities given by
\begin{equation}
  \bry(\bfx) := \brz(\bfx)\brA^{-1}. \quad \bfx\shin\G \label{y:def}
\end{equation}
As linear combinations of $\bfz_i$, the $\bfy_i$ belong to $\bsfH_r$. Thanks to~\eqref{A:zi}, they verify
\begin{equation}
  \alpha_i[\bfy_j]
  = \alpha_i\lsqb \bfz_k \rmA_{kj}^{-1} \rsqb = \rmA_{ik}\rmA_{kj}^{-1} = \delta_{ij}
\end{equation}
(with implicit summation over the repeated subscript $k=1,\ldots,r$), which means
\begin{equation}
  \sfR\bfy_i=\bfuR_i \quad (i=1,\ldots,r), \qquad \text{i.e.}\quad \sfR\bry = \brvR. \label{FR:yS}
\end{equation}
In other words, setting $\bfuS=\bfy_i$ in~\eqref{BC} and enforcing the no-net-force conditions~\eqref{no-net-force} produces on $\G$ a rigid-body velocity equal to the basis function $\bfuR_i$.

The main useful properties of the slip velocities $\bfy_j$ thus defined are gathered in the next proposition, whose proof is given in Appendix~\ref{sec:proof:y:props}:\enlargethispage*{5ex}
\begin{proposition}\label{y:props}
Let the $r$-dimensional subspace $\bsfH_r$ of $\bsfHs$ be defined by $\bsfH_r:=\text{span}(\bfy_1,\ldots,\bfy_r)=\text{span}(\bfz_1,\ldots,\bfz_r)$.
\begin{compactenum}[(i)]
\item The linear mapping $\bfuS\mapsto\bfupal[\bfuS]$ defines a  $\bsfH_r\to\Rbb^r$ bijection. It verifies
\begin{equation}
  \bfuS =\bry\bfal \in \bsfH_r \quad\Longleftrightarrow\quad
  \bfupal[\bfuS] = \bfal\in\Rbb^r \label{H:R}
\end{equation}
for any $\bfal\shin\Rbb^r$, so that any specified rigid-body velocity characterized by $\bfal=[\bfU;\bfOm]$ is generated by setting the slip velocity in~\eqref{BC:swimmer} to $\bfuS=\bry\bfal$.
\item Let $\bfuS=\bry\bfal\in\bsfH_r$ for some $\bfal\in\Rbb^r$. The power loss is then given by $  P(\bfuS) = \bfal\tsup\brA^{-1}\bfal$.
\item For any $(\bfuS,\bfwS)\in\bsfH_r\shtimes\bsfN$, we have $P(\bfuS\shp\bfwS) \geq P(\bfuS)$.
\end{compactenum}
\end{proposition}
In particular, item (iii) implies that the low-dimensional optimization problem~\eqref{JL:opt:lowdim} is equivalent to the original optimization problem~\eqref{JL:opt:generic} and that $P(\bry\bfal)=\bfal\tsup\brA^{-1}\bfal$ is the lower bound of the incurred power loss achieved by all slip-induced motions having the same rigid-body parameters $\bfal\shin\Rbb^r$.

\begin{remark}[link to the minimum dissipation result of~\cite{nasouri:21}]
Fixing $\bfal\shin\Rbb^r$, our lower bound $\bfal\tsup\brA^{-1}\bfal$ of $P$ is equal to that found in~\cite{nasouri:21}, where a different reasoning is followed. To prove this, we note that (a) setting $\bfbe=(\bfA\shp\bfC^{-1})^{-1}\bfal\in\Rbb^r$ gives a flow such that $\alpha_i[\brz\bfbe]=(\bfA\bfbe)_i$, (b) applying the rigid-body velocity $\brvR\brC^{-1}\bfbe$ creates on $\G$ the traction $\bff'=\brfR\bfC^{-1}\bfbe$ whose tangential projection $\bfPi\bff'$ exactly cancels that created by the slip $\brz\bfbe$, see~\eqref{forward:PDE:fD:JL}. The superposition of states (a) and (b) thus induces the rigid-body motion whose parameters are $\brA\bfbe+\brC^{-1}\bfbe=(\bfA\shp\bfC^{-1})(\bfA\shp\bfC^{-1})^{-1}\bfal=\bfal$; it is thus the perfect-slip flow creating $\bfal$. From the no-net-force conditions satisfied by the flow (a), the resistance tensor $\brR_{\text{PS}}$ associated with the perfect-slip flow (a)+(b) for given $\bfal$ is obtained from its no-slip part (b) as $(\brA\shp\brC^{-1})^{-1}$. Recalling Remark~\ref{rem:RNS}, we finally note that
\begin{equation}
  \bfal\tsup\brA^{-1}\bfal = \bfal\tsup\lsqb (\brA+\brC^{-1}) - \brC^{-1} \rsqb^{-1}\bfal
  = \bfal\tsup\lsqb \brR_{\text{PS}}^{-1} - \brR_{\text{NS}}^{-1} \rsqb^{-1}\bfal,
\end{equation}
which, upon adjusting notation, coincides with the bound given by~\cite[eq.~1]{nasouri:21}. Likewise, the efficiency $\eta_{\rmm}$ for fixed $\bfal$ defined by~\cite[eq.~2]{nasouri:21} is given, with the present notations, by $\eta_{\rmm} =  \bfal\tsup ( \brA \shp \brC^{-1} )^{-1} \bfal\, /\, \bfal\tsup \brA^{-1} \bfal$.

In addition to establishing the lower bound of~\cite{nasouri:21} via an alternative route, our treatment provides a computationally constructive approach for evaluating that lower bound and setting slip optimization problems in low-dimension form via the definition and computation of $\brz$ and the resulting matrix $\brA$ and slip velocities $\bry$.
\end{remark}

\subsection{Boundary integral framework}
\label{sec:bie}

We solve the $2r$ flow problems involved in this work, namely the $r$ no-slip Dirichlet BVPs~\eqref{Stokes}-\eqref{BC}
with $\bfuD = \bfuR_\ell$ and the $r$ mixed BVPs~\eqref{forward:PDE:fD:JL}, %---
using an indirect integral equation formulation with the single-layer potential (SLP) {\em ansatz}. Letting $\boldsymbol{G}$ and $\boldsymbol{T}$ denote the Stokeslet and Stresslet fundamental solutions, respectively given by \cite{corona:17,Pozrikidis1992}
\begin{equation}
\boldsymbol{G}(\bfr) =  \frac{1}{8\pi\mu}
       \left(\frac{\bfI}{|\bfr|}
       + \frac{\bfr\otimes\bfr}{|\bfr|^3}\right), \qquad \boldsymbol{T}(\bfr) =  -\frac{3}{4\pi}
\frac{\bfr\otimes\bfr\otimes\bfr}{|\bfr|^5}, \label{stokeslet}
\end{equation}
the velocity field is hence expressed as
\begin{equation}
    \bfu(\bfx) = \mathcal{S}[\bfmu](\bfx) := \iG \boldsymbol{G}(\bfr) \bfmu(\bfy) \dS_y, \qquad \bfr := \bfx - \bfy, \label{eq:slp}
\end{equation}
for an unknown surface density $\bfmu$. The representation~\eqref{eq:slp} satisfies the Stokes equations in $\OO$ and the decay condition at infinity by construction. For the Dirichlet BVPs, enforcing the boundary condition $\bfu = \bfuD$ on $\G$ via the continuity of the SLP across $\G$ yields a Fredholm integral equation of the first kind for $\bfmu$:
\begin{equation}
    \mathcal{S}[\bfmu](\bfx) = \bfuD(\bfx), \qquad \bfx \in \G.
    \label{eq:bie}
\end{equation}
The traction $\bff = -p\bfn + 2\mu\bfD[\bfu]\sip\bfn$ on $\G$ is then recovered from
the density $\bfmu$ using the standard jump relations~\cite{corona:17,Pozrikidis1992} for layer potentials as
\begin{equation}
\bff(\bfx) = -\tfrac{1}{2}\bfmu(\bfx) + \mathcal{K}[\bfmu](\bfx),
\quad\text{where}\quad
\mathcal{K}[\bfmu](\bfx) := \iG \boldsymbol{T}(\bfr) \bfmu(\bfy)\,\mathrm{d}S_y \cdot \bfn(\bfx),
\label{eq:traction}
\end{equation}
and the integral in $\mathcal{K}$ is taken in the principal-value sense.

For mixed BVPs~\eqref{forward:PDE:fD:JL}, the boundary condition prescribes the \emph{tangential traction} and the tangentiality constraint $(\bfv_i - \bfv_i\Rsup)\tsup\bfn = 0$ on $\G$. Starting from the SLP ansatz~\eqref{eq:slp} for the velocity field again, we can derive the following set of boundary integral equations for the unknowns $\{\bfmu_i, \bfv_i \Rsup\}$:
\begin{equation}
    \begin{aligned}
        \bfPi\!\left(-\tfrac{1}{2}\bfmu_i(\bfx) + \mathcal{K}[\bfmu_i](\bfx)\right) &= \bfPi\bff^\alpha_i(\bfx), \\
        \mathcal{S}[\bfmu_i] (\bfx) \cdot \bfn   &= \bfv_i \Rsup \cdot \bfn  ,
    \end{aligned}
    \qquad \bfx \in \G.
    \label{eq:bie_aux}
\end{equation}
The zero net force and torque conditions~\eqref{no-net-force} provide the additional equations necessary to fully determine the unknown density and rigid-body velocities.\enlargethispage*{7ex}

The microswimmer surface $\G$ is parametrized using spherical harmonics of degree $p$, as in
the simulation framework of~\cite{corona:17,shravan:11}: coordinate functions and
surface densities are represented in the basis $\{Y^m_\ell\}_{|\ell|\leq p}$, with
a $(p+1)$-point Gauss--Legendre rule in the polar direction and the trapezoidal rule
in the azimuthal direction providing spectral accuracy for smooth integrands. The weakly singular integral operators $\mathcal{S}$ and $\mathcal{K}$ are numerically evaluated using the fast spherical grid rotation quadrature of~\cite{gimbutas:13}. The resulting dense linear systems are solved iteratively via GMRES.

The scalar products $\lbra \cdot, \cdot \rbra_\G$ entering the assembly of the matrices
$\brC$ and $\brA$ (see~\eqref{CR:expr} and~\eqref{eig:rank:r:JL}) are then evaluated using
the same Gauss--Legendre quadrature.

\subsection{Numerical validations}
\label{sec:valid}

\begin{figure}[b]
    \centering
    \includegraphics[trim=0 30 0 0, clip, width=\textwidth]{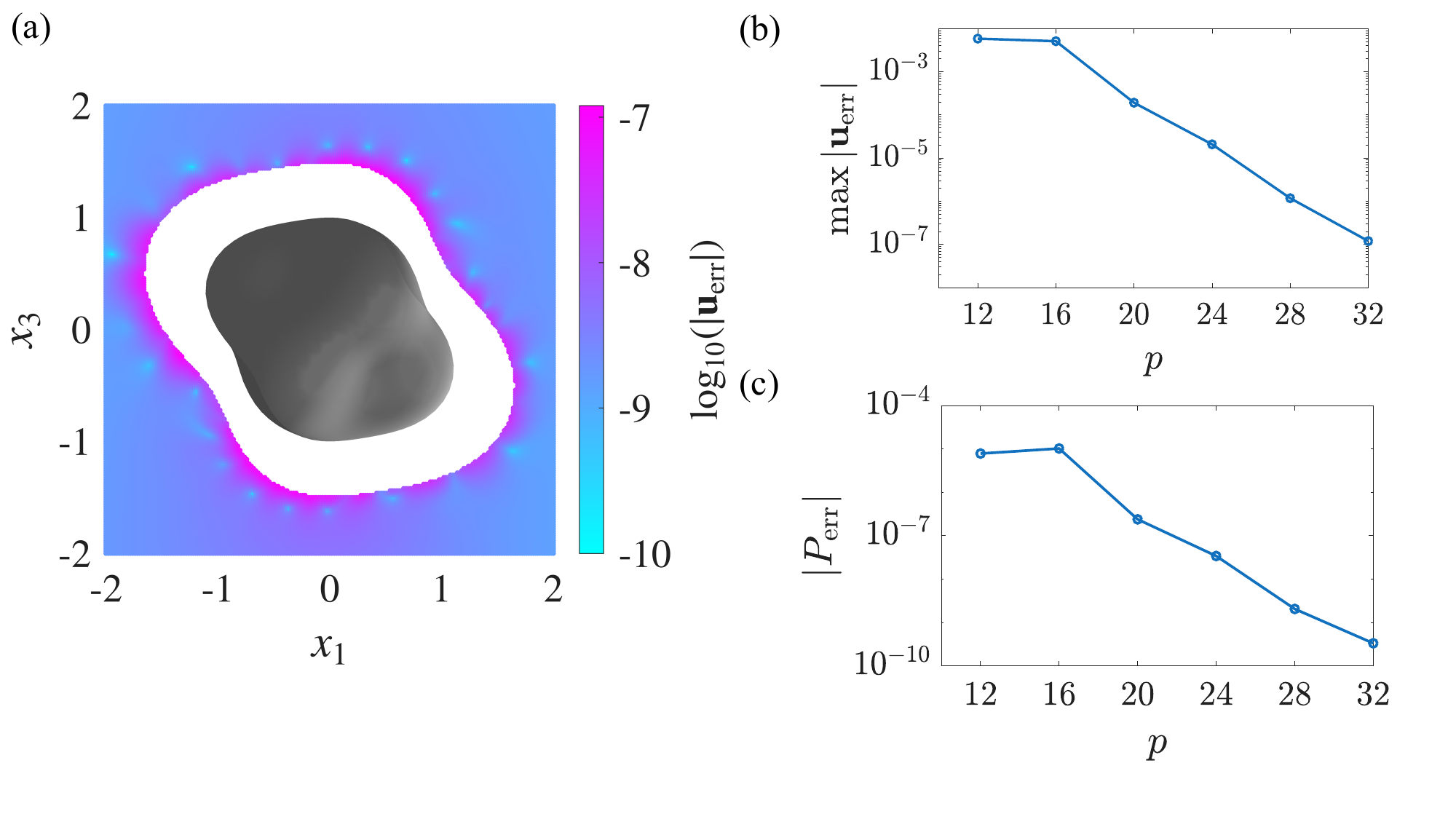}\vspace*{-2ex}
    \caption{Numerical validation of the boundary integral method on a mixed BVP for the swimmer shape defined using spherical coordinates $(r,\theta,\phi)$ by $r(\theta,\phi)=1 + 0.5 \textrm{Re}(Y_4^3(\theta,\phi))$ (with $Y_{\ell}^m$ the complex-valued spherical harmonic of degree $\ell$ and order $m$), the reference flow being that created by a point force $\boldsymbol{F} = [1,1/2,1/3]\tsup$ placed inside the swimmer at $\boldsymbol{x}_0 = [0.1,0.2,-0.3]\tsup$. (a) The absolute error between the exact solution and the numerical solution in the $x_2 = 0$ plane with a total of 1200 Gaussian quadrature points. (b) $L_\infty$-error in the flow field when $2p(p+1)$ Gauss-Legendre points are used to discretize the surface. (c) $L_\infty$-error in the power loss as a function of $p$.}
    \label{fig: convergence}\bigskip

    \includegraphics[trim=0 80 0 30, clip, width=0.8\textwidth]{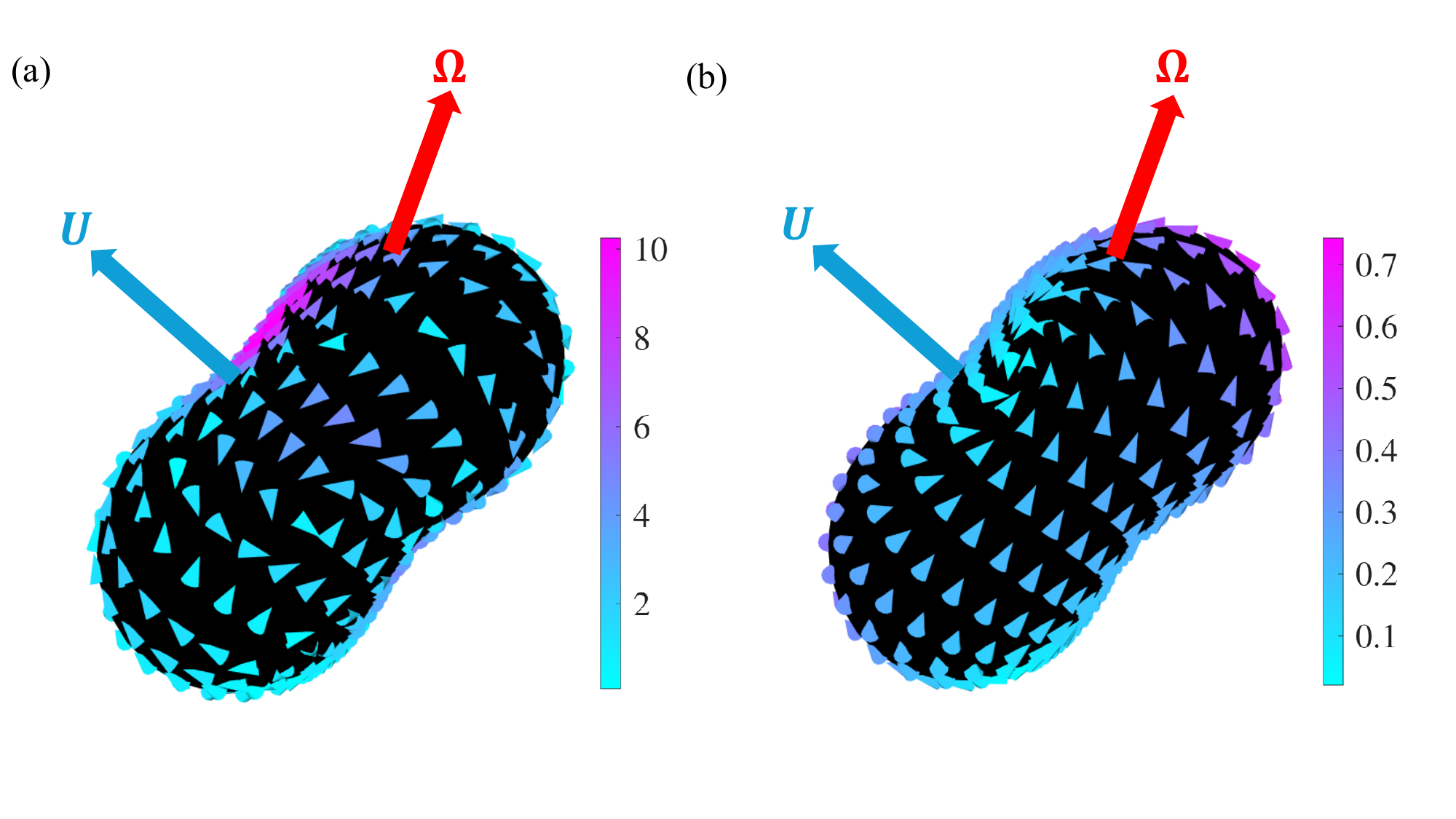}
  \caption{Verification of Proposition \ref{y:props} on a tilted dumbbell whose shape is defined, with notations as in Fig.~\ref{fig: convergence}, by $r(\theta,\phi)=1 + \textrm{Re}(Y_2^1(\theta,\phi)) + 0.1 \textrm{Re}(Y_3^2(\theta,\phi))$. (a) Plot of an arbitrary slip velocity $\bfuS \in \bsfHs$ with corresponding rigid body velocities $\boldsymbol{U} = [0.15,0.74,0.26]\tsup$ and $\boldsymbol{\Omega} = [0.53,0.01,0.92]\tsup$ and power loss $P(\bfuS) = 809.74$. (b) The slip velocity $\bry\bfal \in \bsfH_r$ produces identical rigid body motion $\bfal =[\bfU;\bfOm]$ with reduced power loss $P(\bry\bfal) = 30.95 \leq P(\bfuS).$}\label{fig:lemma-verify}
\end{figure}

We first verify the boundary integral method on a mixed BVP that, like the BVP~\eqref{forward:PDE:fD:JL}, features prescribed normal velocity and tangential traction. An exact solution is generated by a point source $\boldsymbol{F}$ placed at a position $\boldsymbol{x}_0$ within an arbitrarily-shaped swimmer, its flow and traction fields being given by 
\begin{equation}
\boldsymbol{u}_{\mathrm{exa}} = \boldsymbol{G}(\boldsymbol{x}-\boldsymbol{x}_0) \boldsymbol{F}, \quad \boldsymbol{f}[\boldsymbol{u}_{\mathrm{exa}}] = \boldsymbol{T}(\boldsymbol{x} - \boldsymbol{x}_0) \boldsymbol{F} \cdot \boldsymbol{n}(\boldsymbol{x}), \label{sol:exact}
\end{equation}
where $\boldsymbol{G}$ and $\boldsymbol{T}$ are again the Stokeslet and Stresslet~\eqref{stokeslet}. We then prescribe $\bfn\tsup\boldsymbol{u}_{\mathrm{exa}}$ and $\bfPi \boldsymbol{f}[\boldsymbol{u}_{\mathrm{exa}}]$ on $\G$ and solve numerically the resulting mixed BVP using the boundary integral method described in Section~\ref{sec:bie}. Numerical results for the flow fields and the power loss are seen in Fig.~\ref{fig: convergence} to converge to the exact solution~\eqref{sol:exact} as the quadrature parameter $p$ is increased.

Then, Figure~\ref{fig:lemma-verify} presents a numerical verification of Proposition~\ref{y:props}, where an arbitrarily-chosen slip velocity profile $\bfuS$ is first applied on a different arbitrarily-shaped swimmer, inducing a power loss $P(\bfuS) = 809.74$. The resulting rigid-body motion parameters $\bfal=[\bfU,\bfOm]$ are then obtained using~\eqref{R:expr}. Next, the slip velocity $\bfuS=\bry\bfal$ is applied to the swimmer. The rigid-body parameters evaluated from that solution are found to coincide, as predicted, with $\bfal$ (i.e. $\bfupal[\bfuS]=\bfal$), while a substantially lower power loss $P(\bry\bfal)= 30.95$ is achieved.

\section{Low-dimensional power loss minimization problem}
\label{sec:solv:opt}

\subsection{Generic optimization problem}
\label{JW:opt:generic}

Recalling Proposition~\ref{y:props}, any $\bfuS\shin\bsfH_r$ may be written in terms of the translation and rotation velocities $\bfU,\,\bfOm$ arising from enforcing the no-net-force conditions as
\begin{equation}
  \bfuS = \bry\bfal = \bryU\bfU + \bryO\bfOm, \qquad
  \bfal = [\bfU;\bfOm] \in\Rbb^r. \label{aux15}
\end{equation}
The incurred power loss is then given (see item (ii) of Prop.~\ref{y:props}) by
\begin{equation}
  P(\bfuS)
 = \bfal\tsup\brA^{-1}\bfal
 = \bfU\tsup\brZUU\bfU + 2\bfU\tsup\brZUO\bfOm + \bfOm\tsup\brZOO\bfOm, \qquad \brA^{-1}
 = \begin{bmatrix} \brZUU & \brZUO \\ \brZOU & \brZOO \end{bmatrix}
\end{equation}
(with $\brZOU=\brZUO\tsup$) using the above 
partition of $\brA^{-1}$ induced by the partition $\bfal = [\bfU;\bfOm]$. Introducing the parametrization~\eqref{rigid:param}, the power loss is then expressed in terms of $s,\bfV,\bfW$ by
\begin{align}\label{aux20}
  \hspace*{2em}
  P(\bfuS)
 &= \Pcal(s,\bfV,\bfW) \\
 &= \bfV\Tsup\brZUU\bfV + 2\bfV\Tsup\lpar s\brZUO + \brZUU \rpar\bfW
 + \bfW\Tsup\lpar s^2\brZOO + 2s\brZUO + \brZUU \rpar \bfW,
\end{align}
and problem~\eqref{JL:opt:lowdim} takes the form
\begin{equation}
  \min_{(s,\bfV,\bfW)} \Pcal(s,\bfV,\bfW) \qquad \text{subject to} \qquad
  \bfV\tsup\bfW=0,\ |\bfW|=1. \label{JL:opt:rdim}
\end{equation}
With reference to the discussion of parametrization~\eqref{rigid:param}, $\Pcal$ is a well-behaved function of $(s,\bfV,\bfW)$ that remains defined at values $(s,\bfV,\bfW)$ that do not produce a valid net motion, see Sec.~\ref{sec:AUO=0}.\enlargethispage*{-1ex}

It is in fact convenient to treat problem~\eqref{JL:opt:rdim} as two nested minimizations, where $\bfW=\bfW^{\star}$ solves the reduced problem
\begin{equation}
  \min_{\bfW\in\Sigma} \PcalH(\bfW), \qquad \PcalH(\bfW) := \Pcal(s[\bfW],\bfV[\bfW],\bfW)
\label{JW:min:angular}
\end{equation}
($\Sigma$ being the unit sphere) and $s[\bfW],\bfV[\bfW]$ solve, given $\bfW$, the partial minimization problem
\begin{equation}
  \min_{(s,\bfV)} \Pcal(s,\bfV,\bfW) \qquad \text{subject to} \qquad \bfV\tsup\bfW=0. \label{JL:opt:rdim:part}
\end{equation}
Recalling parametrization~\eqref{rigid:param} and Proposition~\ref{y:props}, the corresponding optimal slip velocity is finally obtained in terms of $(s,\bfV,\bfW)$ solving problem~\eqref{JL:opt:rdim} as
\begin{equation}
  \bfuS = \bryU(\bfW\shp\bfV) + s\bryO\bfW. \label{uS:opt}
\end{equation}

\subsection{Partial minimization with $\bfW$ fixed}
\label{generic:min}

Problem~\eqref{JL:opt:rdim:part} involves a quadratic objective function $(s,\bfV)\mapsto \Pcal(s,\bfV,\bfW)$ and a linear constraint, and thus has a closed-form solution given in the following lemma (whose proof is given in Appendix~\ref{sec:proof:partial:min}):
\begin{lemma}\label{partial:min}
Let $\bfW\shin\Rbb^3$ be fixed. The solution $(s,\bfV)$ of problem~\eqref{JL:opt:rdim:part} is given by
\begin{equation}
% \MoveEqLeft[10]{ 
\label{s:mu:expr}
  s = -\inv{D}|\bfW|^2\AUO, \quad \bfV = \Lpar \inv{D}|\bfW|^2\AOO \brZUU^{-1}-s\brZUU^{-1}\brZUO-\bfI \Rpar \bfW,
\end{equation}
with the $\bfW$-dependent scalars $\AUU,\,\AUO,\,\AOO$ and $D$ given by
\begin{equation}
\label{aux21}
\begin{aligned}
  \AUU &:= \bfW\tsup\brZUU^{-1}\bfW, &\qquad
  \AUO &:= \bfW\tsup\brZUU^{-1}\brZUO\bfW, \\
  \AOO &:= \bfW\tsup \lpar \brZOO - \brZOU\brZUU^{-1}\brZUO \rpar\bfW, &
  D &:= \AUU\AOO + \AUO^2
\end{aligned}
\end{equation}
(noting that $\AUU>0$ and $\AOO>0$ by virtue of $\brA$ being symmetric positive definite). The corresponding rigid-body parameters and power loss are given by
\begin{equation}
  \bfU = \bfV + \bfW, \qquad
  \bfOm = s\bfW, \qquad
  \PcalH(\bfW) = \inv{D}|\bfW|^4\AOO. \label{JW:min:part}
\end{equation}
\end{lemma}
The foregoing partial minimization is sufficient if the net motion direction $\bfW$ is given \emph{a priori}. We also note that $\PcalH(\bfW)$ is homogeneous with degree 2 in $\bfW$, so that $\Sigma$ may be reduced to the half-sphere in problem~\eqref{JW:min:angular}. On the other hand, solving problem~\eqref{JW:min:angular} is useful e.g. for finding the optimal orientation of the swimmer body relative to a given net motion direction $\bfW$.

\subsection{Particular case $\AUO=0$}
\label{sec:AUO=0}

The case where the chosen net direction $\bfW$ satisfies $\AUO(\bfW)=0$ (see~\eqref{aux21}) warrants additional discussion. Applying Lemma~\ref{partial:min} to this case (and assuming $|\bfW|=1$ here without detriment) gives
\begin{equation}
  s=0, \quad \bfOm=\bfze, \quad \bfU= \brZUU^{-1}\bfW/\AUU, \quad \PcalH(\bfW)=1/\AUU. \label{s=0}
\end{equation}
In the exceptional event where $\bfW$ is an eigenvector of $\brZUU^{-1}$, i.e. $\brZUU^{-1}\bfW=\eta\bfW$ for some $\eta>0$, we have $\AUU=\eta$ and hence $\bfU=\bfW$, and a physically-consistent translation motion satisfying~(\ref{W:def}b) is obtained. Otherwise, the above solution~\eqref{s=0} is inconsistent with~\eqref{W:def}.

To interpret the latter situation, we observe that freezing $s$ in the partial minimization yields $\Pcal=1/\AUU + s^2\AOO>\PcalH(\bfW)$ (as found by adapting the proof of Lemma~\ref{partial:min}). The inconsistent situation ($\AUO=0$, $\bfW$ not an eigenvector of $\brZUU^{-1}$) thus is a limiting case where, as $s\to0$, $\bfW$ results from helical centroid trajectories whose radius $|\bfV|/|s|$ diverges, while the expended power decreases towards the lower bound $\PcalH(\bfW)$. Similarly, for $\bfW$ chosen such that $\AUO(\bfW)$ is small (but nonzero), Lemma~\ref{partial:min} yields a small value of $s$ while $\bfV=(\frac{|\bfW|^2}{D}\brZUU^{-1}-\brI)\bfW+O(s)$. The radius $|\bfV|/|s|$ of that optimal helical motion again diverges as $\AUO\to0$. Such small-$s$ or small-$\AUO$ situations are certainly undesirable from a physical standpoint.

\subsection{Reduced minimization over $\bfW$}

Problem~\eqref{JW:min:angular} does not appear to be solvable in closed form, as $D(\bfW)$ is not quadratic in $\bfW$. An iterative numerical algorithm can instead be applied, with $\PcalH(\bfW)$  evaluated for each trial value of $\bfW$ using Lemma~\ref{partial:min}. Such process may encounter intermediate values of $\bfW$ for which Lemma~\ref{partial:min} yields non-physical values $s\sheq0,\,\bfV\shneq\bfze$. However, problem~\eqref{JW:min:angular} verifies the following lemma (proved in Section~\ref{AUO:opt:proof}), which ensures that any solution $(s,\bfV,\bfW)$ to the complete minimization~\eqref{JL:opt:rdim} must produce a valid net motion~\eqref{W:def}, i.e. the undesirable situation $s\sheq0,\,\bfV\shneq\bfze$ cannot occur at the global optimum:
\begin{lemma}\label{AUO:opt}
Let $\bfW^{\star}\shin\Sigma$ verify the first-order optimality conditions for problem~\eqref{JW:min:angular}. If $\bfW^{\star}$ is also such that $\AUO(\bfW^{\star})=0$, then it is an eigenvector of $\brZUU^{-1}$.
\end{lemma}

The following question naturally arises: can optimal motions resulting from the minimization~\eqref{JL:opt:rdim} be helical (and thus feature non-zero rotations)? The forthcoming discussion of Sec.~\ref{CA:sym} shows that helical globally optimal motions are precluded for swimmers with high enough symmetry (in particular axisymmetry), and otherwise may arise only when $\bfW^{\star}$ solving~\eqref{JL:opt:rdim} does not lie in a plane of symmetry. Numerical evidence shows that helical optimal motions can indeed occur, as demonstrated next.\medskip

\paragraph{Numerical example: helical globally-optimal motion}

\begin{figure}[t]
    \centering
    \includegraphics[trim=0 0 180 0, clip, width=\textwidth]{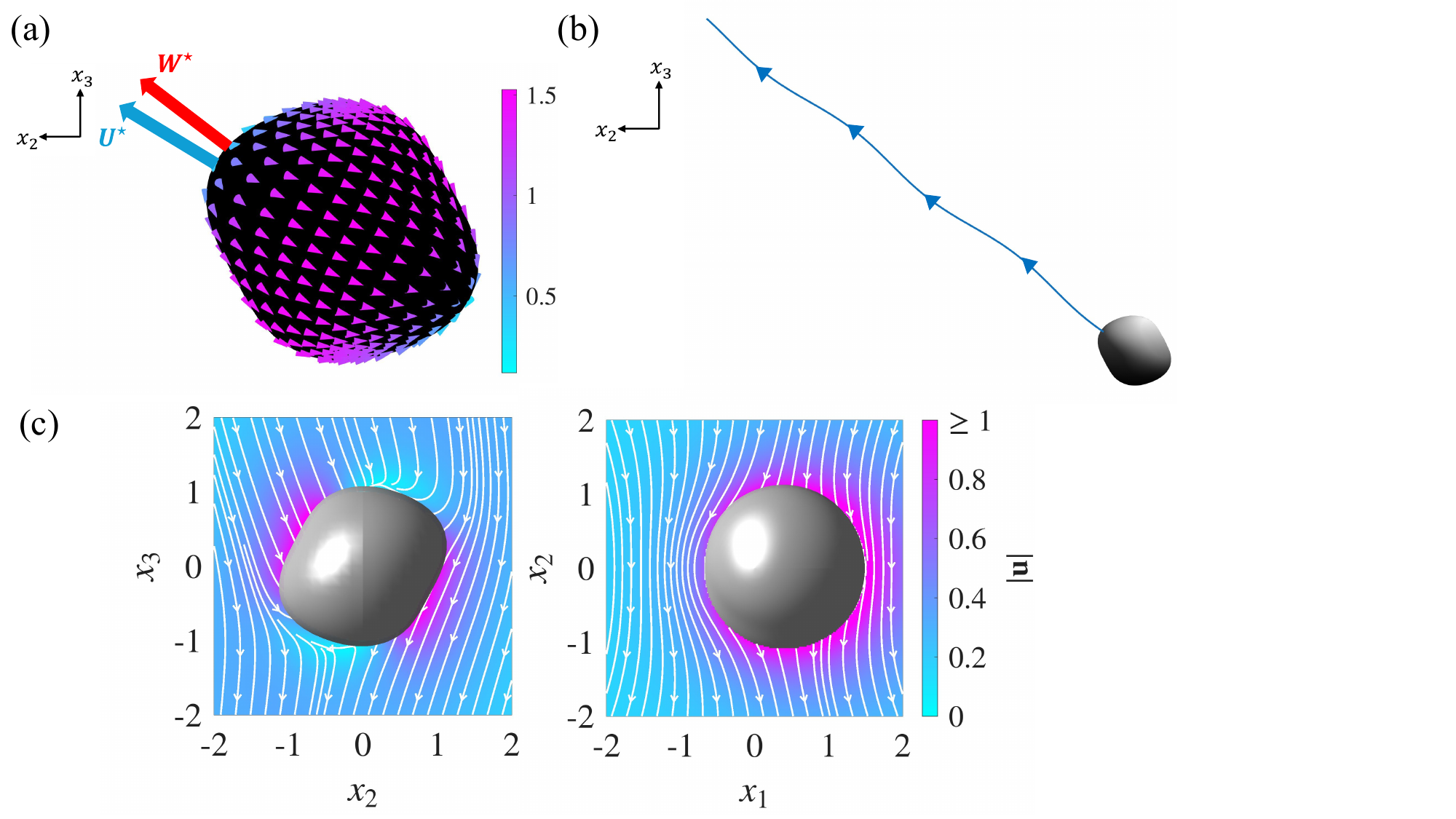}
    \caption{Global optimization of the shape defined by $r(\theta,\phi)=\exp(0.4(\sin\theta\cos\phi + \sin^4\theta \cos\theta\sin\phi))$, with notations as in Fig.~\ref{fig: convergence}. (a) Optimal slip velocity field plotted (with color bar representing $|\boldsymbol{v}^\textrm{S}|$) together with the optimal net motion direction $\boldsymbol{W}^{\star} = -\boldsymbol{\Omega^{\star}}/|\boldsymbol{\Omega^{\star}}| =  (0,0.79,0.61)  $ (red arrow) and the direction of the translational velocity $\boldsymbol{U^{\star}} = (0,0.88,0.50)$ (blue arrow). The nonzero angular velocity $\boldsymbol{\Omega^{\star}} = (0,-0.26, -0.20)$ results in globally optimal helical motion.  (b) Resulting helical motion (with swimmer enlarged for illustrative purposes). (c) Body frame flow field streamlines in the $x_1 = 0$ (left) and $x_3 = 0$ (right) planes.  The optimization resulted in a global minimum power loss of $\widehat{\mathcal{P}}(\boldsymbol{W}^\star) = 35.54 $.}
    \label{fig: skewboomerang_globalopt}
\end{figure}

The global optimization is applied to the shape defined by $r(\theta,\phi)=\exp(0.4(\sin\theta\cos\phi + \sin^4\theta \cos\theta\sin\phi))$, where $(r,\theta,\phi)$ are spherical coordinates, with all flow fields and resulting matrices computed as described in Section~\ref{sec:bie}. This yields the optimal direction of net motion $\bfW^{\star}=(0,0.79,0.61)$, the translational velocity $\boldsymbol{U^{\star}} = (0,0.88,0.50)$ and the angular velocity $\bfOm^{\star} = (0,-0.26, -0.20)$; in particular, $\bfU^{\star}$ and $\bfW^{\star}$ are not collinear, see Figure~\ref{fig: skewboomerang_globalopt}. The obtained global minimum power loss is $\widehat{\mathcal{P}}(\bfW^\star) = 35.54 $.
\enlargethispage*{5ex}

\subsection{Spinning or straight rigid-body motions}
\label{sec:special}

The power loss minimization~\eqref{JL:opt:generic} may be restricted to spinning or straight rigid-body motions, i.e. to case (b) of~\eqref{W:def}. This simpler problem is addressed by setting $\bfV=\bfze$ and ignoring the (now moot) orthogonality constraint in the parametrization~\eqref{rigid:param}, and seeking $s\shin\Rbb$ that minimizes $s\mapsto\Pcal(s,\bfze,\bfW)$ with $\Pcal$ still given by~\eqref{aux20}. This results in
\begin{equation}
  s^{(b)} =  - \frac{\BUO}{\BOO}, \qquad
  \bfuS = \lpar \bryU + s^{(b)}\bryO \rpar \bfW, \qquad
  P^{(b)} = \BUU - \frac{\BUO^2}{\BOO} \label{JW:O}
\end{equation}
with $\BUU := \bfW\tsup \brZUU \bfW$, $\BUO := \bfW\tsup \brZUO \bfW$ and $\BOO := \bfW\tsup \brZOO \bfW$. If $\BUO\shneq0$, the above solution produces a spinning rigid-body motion ($\bfOm\shneq\bfze$). Since this constitutes a partial minimization for problem~\eqref{JL:opt:rdim:part}, we must have $P^{(b)} \geq \PcalH(\bfW)$. When $\BUO\sheq0$, a straight rigid-body motion is obtained.\enlargethispage*{1ex}

If only straight (i.e. rotationless) rigid-body motions are considered, the net motion direction with least power loss solves the Rayleigh quotient minimization
\begin{equation}
  \min_{\bfW\in\Rbb^3,\,|\bfU|=1} \BUU(\bfW)
\end{equation}
and is thus the eigenvector $\bfW^{\star}$ associated with the smallest eigenvalue $\eta^{\star}>0$ of the eigenvalue problem $\brZUU\bfU-\eta\bfU=\bfze$. The optimal slip velocity $\bfuS{}^{\star}$ and power loss $P^{\star}$ are hence given by
\begin{equation}
  \bfuS{}^{\star} = \bryU\bfW^{\star}, \qquad P^{\star} = P(\bfuS{}^{\star}) = \eta^{\star}.
\end{equation}

\subsection{Slip optimization algorithm}
\label{sec:algo:generic}

We obtain the following solution procedure for solving problem~\eqref{JL:opt:generic}, which requires $2r$ flow solutions.\medskip
\begin{compactenum}[(i)]
\item Solve the $r$ problems~\eqref{Stokes}-\eqref{BC} with $\bfuD=\bfuR_{\ell}$ ($1\shleq \ell\shleq r$), obtain tractions $\bff\Rsup_{\ell}$.
\item Set up matrix $\bfCr\in\Rbb^{r\times r}\sym$ given in~\eqref{CR:expr} and rigid body extraction tractions $\bff_i^{\alpha}$ given by~\eqref{R:expr}.
\item Solve the $r$ problems~\eqref{forward:PDE:fD:JL}, obtain $\bfz_i$ given by~\eqref{zi:def:JL} and tractions $\bff_i$. Set up $\brz(\bfx)=[\bfz_1(\bfx),\ldots,\bfz_r(\bfx)]$.
\item Set up matrix $\brA\in\Rbb^{r\times r}\sym$ given by $A_{ij}=\lbra \bff_i,\bfz_j \rbra_{\G}$.
\item Compute $\brZ=\brA^{-1}$ and $\bry(\bfx)=\brz(\bfx)\brZ$.
\item Solve the reduced minimization problem~\eqref{JW:min:angular} by iteration on $\bfW\shin\Sigma$. For each trial value of $\bfW$ encountered, apply Lemma~\ref{partial:min} and evaluate $\PcalH(\bfW)$ using~\eqref{JW:min:part}.
\item Upon convergence, evaluate the optimal slip velocity $\bfuS\shin\bsfH_r$ using~\eqref{aux15} with $\bfU,\bfOm$ given by~\eqref{JW:min:part}.\medskip
\end{compactenum}
Steps (i) to (iv) are performed using the boundary integral framework described in Section~\ref{sec:bie}. The remaining steps (v) to (vii) are negligible in cost relative to the flow solves as they involve only $r \shtimes r$ dense linear algebra and nonlinear methods on $r$-vectors. To solve only the partial minimization~\eqref{JL:opt:rdim} given $\bfW$, perform steps (i) to (v) and apply Lemma~\ref{partial:min}, taking into account the discussion of Sec.~\ref{sec:AUO=0} if $\AUO(\bfW)=0$.

\section{Structure of $\brC,\brA$ for swimmers with symmetry}
\label{CA:sym}

We address the cases for which $\G$ has one, two or three planes of symmetry, and also that of dihedral symmetry (which in particular includes all axisymmetric swimmers); those cases are thereafter labeled S1, S2, S3, D and (for axisymmetry) A. Their outcome mainly consists of the patterns obeyed by matrices $\brC,\brA$ due to those symmetries, whose consequences on the result of (partial or complete) slip optimization are as follows:
\begin{proposition}\label{prop:optim:symm}
Let $\bfW^{\star}$ solve~\eqref{JL:opt:rdim}, i.e. define a globally optimal motion.

1. In cases A and S3, $\bfW^{\star}$ must lie in a plane of symmetry, and the corresponding globally optimal motion is rotationless.

2. In the other cases, if $\bfW$ lies in a symmetry plane, the partial minimization~\eqref{JL:opt:rdim:part} yields a rotationless net motion whenever its solution is consistent (see Sec.~\ref{sec:AUO=0}). In particular, if $\bfW^{\star}$ lies in a symmetry plane, the globally optimal motion is rotationless.
\end{proposition}
In what follows, we justify Prop.~\ref{prop:optim:symm} by relying on the following premise: any relevant boundary-value problem with datum $\bfd$ (e.g. $\bfd=\bfuS$ for problem~\eqref{Stokes}-\eqref{BC}) verifies
\begin{equation}
    \bfd(\brs\dotp)=\brs\bfd(\dotp) \ \text{on } \G \quad \implies \quad\bfu(\brs\dotp)=\brs\bfu(\dotp), \ \  \bff(\brs\dotp)=\brs\bff(\dotp) \ \text{on } \G, \label{BVP:sym}
\end{equation}
where $\brs$ is an isometry of $\Rbb^d$ under which $\OO$ is invariant ($\bfx\shin\OO\,\Longleftrightarrow\,\brs\bfx\shin\OO$). Implication~\eqref{BVP:sym} can be proved using the general and rigorous methods for exploiting geometrical symmetry in linear PDEs developed in e.g.~\cite{all:98,bos:86}; we omit here this step to avoid overly lengthy developments.

\subsection*{One symmetry plane (case S1)} Assume for definiteness that the swimmer is invariant under the plane symmetry $\brs_3$ w.r.t. the $x_3=0$ plane (i.e. $\bfx=(x_1,x_2,x_3)\shin\OO \,\Longleftrightarrow\, \brs_3\bfx:=(x_1,x_2,-x_3)\shin\OO$). The behavior under symmetry of the rigid-body basis functions is then easily checked to be given by
\begin{equation}
  \bfuR_{\ell}(\brs_3\bfx)=\brs_3\bfuR_{\ell}(\bfx) \quad (\ell=1,2,6); \qquad
  \bfuR_{\ell}(\brs_3\bfx)=-\brs_3\bfuR_{\ell}(\bfx) \quad (\ell=3,4,5). \label{sym:S3}
\end{equation}
Invoking~\eqref{BVP:sym}, the traction fields $\bffR_{\ell}$ solving problem~\eqref{Stokes}-\eqref{BC} with $\bfuD=\bfuR_{\ell}$ also behave according to~\eqref{sym:S3} under the action of $\brs_3$.

Turning to the evaluation of $\brC$, we let $G$ be a subset of $\G$ such that $G\shcap \brs_3(G)=\emptyset$, $\overline{G\shcup\brs_3(G)}=\G$, which allows to write
\begin{equation}
  \lbra f,g \rbra_{\G}
 = \lbra f,g \rbra_{G} + \lbra f,g \rbra_{\brs_3(G)}
 = \lbra f,g \rbra_{G} + \lbra f(\brs_3\dotp),g(\brs_3\dotp) \rbra_{G} \label{int:split}
\end{equation}
(having set $\bfx=\brs_3\bfx'$ on $\brs_3(G)$ and used that $\brs_3$ is an isometry). Applying the above rule to the entries of $\brC$ defined by~\eqref{CR:expr} and exploiting~\eqref{sym:S3}, we obtain e.g.
\begin{equation}\label{Ckl:sym}
\begin{aligned}
  \rmC_{11}
 &= \lbra \bfuR_1, \bffR_1 \rbra_{G} + \lbra \brs_3\bfuR_1, \brs_3\bffR_1 \rbra_{G}
 &&= 2\lbra \bfuR_1, \bffR_1 \rbra_{G} \\
  \rmC_{31}
 &= \lbra \bfuR_3, \bffR_1 \rbra_{G} + \lbra \brs_3\bfuR_3, \brs_3\bffR_1 \rbra_{G}
 &&= 0,
\end{aligned}
\end{equation}
with all other entries of $\brC$ treated similarly, and find ($\brC$ being symmetric) that the following entries vanish:
\begin{equation}
  \rmC_{k\ell} = \rmC_{\ell k} = 0 \quad (k=1,2,6;\ \ell=3,4,5).
\end{equation}
This makes $\brC$ block-diagonal up to row and column permutations, so that corresponding entries of $\brC^{-1}$ also cancel:
\begin{equation}
  \rmC^{-1}_{k\ell} = \rmC^{-1}_{\ell k} = 0 \quad (k=1,2,6;\ \ell=3,4,5) \label{S:S3}
\end{equation}

Next, thanks to the cancellations~\eqref{S:S3} and the relations~\eqref{sym:S3} also obeyed by the tractions $\bffR_{\ell}$, the rigid body extraction tractions $\bff_i^{\alpha}$ given by~\eqref{R:expr} are found to also obey the relations~\eqref{sym:S3}. For example, we have
\begin{equation}
  \bff_3^{\alpha} = \bffR_3\brC^{-1}_{33} + \bffR_4\brC^{-1}_{43} + \bffR_5\brC^{-1}_{53} \quad
  \implies \quad \bff_3^{\alpha}(\brs_3\dotp)=-\brs_3\bff_3^{\alpha}(\dotp).
\end{equation}
The resulting behavior under $\brs_3$ of the tangential tractions $\bfPi\bff^{\alpha}_i$ acting as loadings in problems~\eqref{forward:PDE:fD:JL} in turn allow to infer by~\eqref{BVP:sym} that their solutions $\bfv_i,\bff[\bfv_i]$ also obey~\eqref{sym:S3}. Hence, evaluation patterns such as~\eqref{Ckl:sym} apply in the same way to $\brA$, whose entries together with those of $\brA^{-1}=\brZ$ thus also verify~\eqref{S:S3}.

The outcome of this analysis is that the matrices $\brC,\brA$ and their inverses all have the structure
\begin{equation}
  \brC,\brA,\brC^{-1},\brZ
  = \begin{bmatrix}
  \times & \times & 0 & 0 & 0 & \times \\
  \times & \times & 0 & 0 & 0 & \times \\
  0 & 0 & \times & \times & \times & 0 \\
  0 & 0 & \times & \times & \times & 0 \\
  0 & 0 & \times & \times & \times & 0 \\
  \times & \times & 0 & 0 & 0 & \times \end{bmatrix} \label{1plane}
\end{equation}
(where $\times$ indicates entries that do not vanish because of geometrical symmetry). The pattern~\eqref{1plane}  implies in particular that $\AUO(\bfW)=0$ for any $\bfW=(W_1,W_2,0)$ lying in the plane of symmetry. However, $\bfW$ is not necessarily an eigenvector of $\brZUU$, so may not define a consistent partial optimum in the sense of Sec.~\ref{sec:AUO=0}.\enlargethispage*{1ex}

\subsection*{Two orthogonal symmetry planes (S2)} Assume for definiteness that the swimmer is invariant under the plane symmetries $\brs_1,\brs_2$ w.r.t. $x_1=0$ and $x_2=0$ planes (i.e. $\bfx=(x_1,x_2,x_3)\shin\OO \,\Longleftrightarrow\, \brs_1\bfx:=(-x_1,x_2,x_3)\shin\OO,\ \text{and}\ \brs_2\bfx:=(x_1,-x_2,x_3)\shin\OO$). The behavior under symmetry of the rigid-body basis functions is then easily checked to be given (with $\brt:=\brs_1\brs_2$ denoting the half-turn about $\bfe_3$) by
\begin{equation}
\begin{aligned}
  \bfuR_{\ell}(\brs_1\bfx) &=\brs_1\bfuR_{\ell}(\bfx) \quad (\ell=2,3,4), &\qquad
  \bfuR_{\ell}(\brs_1\bfx) &=-\brs_1\bfuR_{\ell}(\bfx) \quad (\ell=1,5,6), \\
  \bfuR_{\ell}(\brs_2\bfx) &=\brs_2\bfuR_{\ell}(\bfx) \quad (\ell=1,3,5), &\qquad
  \bfuR_{\ell}(\brs_2\bfx) &=-\brs_2\bfuR_{\ell}(\bfx) \quad (\ell=2,4,6), \\
  \bfuR_{\ell}(\brt\bfx) &=\brt\bfuR_{\ell}(\bfx) \quad (\ell=3,6), &\qquad
  \bfuR_{\ell}(\brt\bfx) &=-\brt\bfuR_{\ell}(\bfx) \quad (\ell=1,2,4,5),
\end{aligned} \label{sym:S12}
\end{equation}
Reasoning along the same lines as in case S1, we obtain the structure of $\brC,\brC^{-1}$, then show that the $\bff^{\alpha}_i$ behave as~\eqref{sym:S12} under symmetry, then infer that $\bfv_i,\bff[\bfv_i]$ behave likewise, to finally obtain the structure of $\brA,\brZ$. The added plane symmetry results in additional vanishing entries in all relevant matrices, whose pattern is
\begin{equation}
  \brC,\brA,\brC^{-1},\brZ
  = \begin{bmatrix}
  \times & 0 & 0 & 0 & \times & 0 \\
  0 & \times & 0 & \times & 0 & 0 \\
  0 & 0 & \times & 0 & 0 & 0 \\
  0 & \times & 0 & \times & 0 & 0 \\
  \times & 0 & 0 & 0 & \times & 0 \\
  0 & 0 & 0 & 0 & 0 & \times \end{bmatrix} \label{2plane}
\end{equation}
This result is arrived at with the help of choosing a subset $G$ of $\G$ such that $G\shcap \brs_1(G)\shcap \brs_2(G)\shcap (\brt(G))=\emptyset$, $\overline{G\shcup\brs_1(G)\shcup\brs_2(G)\shcup(\brt(G))}=\G$ and writing all scalar-product integrals as
\begin{equation}
  \lbra f,g \rbra_{\G}
 = \lbra f,g \rbra_{G} + \lbra f(\brs_1\dotp),g(\brs_1\dotp) \rbra_{G}
   + \lbra f(\brs_2\dotp),g(\brs_2\dotp) \rbra_{G}
   + \lbra f(\brt\dotp),g(\brt\dotp) \rbra_{G}. \label{int:split:S12}
\end{equation}
The pattern~\eqref{2plane} implies in particular that $\AUO(\bfW)=0$ for $\bfW=(W_1,0,W_3)$ or $\bfW=(0,W_2,W_3)$, i.e. for any direction $\bfW$ lying in either symmetry plane. As in case S1, $\bfW$ is not necessarily an eigenvector of $\brZUU$, i.e. may not define a consistent partial optimum.

\subsection*{Dihedral symmetry, rotational symmetry (cases D, A)} Assume that $\G$ is invariant under the plane symmetries $\brs_1,\brs_2$ (defined as above), and also under the rotation $\brr$ of angle $\pi/2$ about $\brO\bfe_3$ (making $\G$ invariant under the $D_4$ dihedral group -- note that $\brs_1\brs_2=\brt=\brr^2$, making the foregoing prescription redundant due to $\brs_2=\brs_1\brr^2=\brs_1\brt$).  Then, the pattern~\eqref{2plane} still holds, with some additional links between nonzero entries. In fact, due to the additional geometrical invariance of $\G$ under $\brr$, we have
\begin{align}
  \MoveEqLeft[5]{\bfuR_2 = \brr\bfuR_1(\brr^{-1}\dotp), \quad \bfuR_5 = \brr\bfuR_4(\brr^{-1}\dotp) } \notag \\ 
  &\implies \quad
  \bffR_2 = \brr\bffR_1(\brr^{-1}\dotp), \quad \bffR_5 = \brr\bffR_4(\brr^{-1}\dotp),
\end{align}
which in turn, by virtue of $\brr$ being unitary, provides
\begin{equation}
\begin{aligned}
  C_{22} &= \lbra \bffR_2,\bfuR_2 \rbra_{\G}
 = \lbra \brr\bffR_1(\brr^{-1}\dotp),\brr\bfuR_1(\brr^{-1}\dotp) \rbra_{\G} &&= C_{11}, \\
  C_{55} &= \lbra \bffR_5,\bfuR_5 \rbra_{\G}
 = \lbra \brr\bffR_4(\brr^{-1}\dotp),\brr\bfuR_4(\brr^{-1}\dotp) \rbra_{\G} &&= C_{44}
\end{aligned}
\end{equation}
and
\begin{align}
\MoveEqLeft[10]{
  C_{24} = \lbra \bffR_2,\bfuR_4 \rbra_{\G}
  = \lbra \brr\bffR_1(\brr^{-1}\dotp),\bfuR_4 \rbra_{\G}
  = \lbra \brt\bffR_1(\brt^{-1}\dotp),\brr\bfuR_4(\brr^{-1}\dotp) \rbra_{\G} } \notag \\
  &= -\lbra \bffR_1,\brr\bfuR_4(\brr^{-1}\dotp) \rbra_{\G} = -\lbra \bffR_1,\bfuR_5 \rbra_{\G} = -C_{15}
\end{align}
so that the matrices $\brC,\brA$ and $\brC^{-1},\brZ$ have the following structure:
\begin{align}
\MoveEqLeft[5]{
\label{axisym}
  \qquad \brC,\brA
  = \begin{bmatrix*}[r]
  a & 0 & 0 & 0 & b & 0 \\
  0 & a & 0 & -b & 0 & 0 \\
  0 & 0 & c & 0 & 0 & 0 \\
  0 & -b & 0 & e & 0 & 0 \\
  b & 0 & 0 & 0 & e & 0 \\
  0 & 0 & 0 & 0 & 0 & d \end{bmatrix*} \quad\implies\quad
    \brC^{-1},\brZ
  = \begin{bmatrix*}[r]
  a'\!\! & 0 & 0 & 0 & b'\!\! & 0 \\
  0 & a'\!\! & 0 & -b'\!\! & 0 & 0 \\
  0 & 0 & c'\!\! & 0 & 0 & 0 \\
  0 & -b'\!\! & 0 & e'\!\! & 0 & 0 \\
  b'\!\! & 0 & 0 & 0 & e'\!\! & 0 \\
  0 & 0 & 0 & 0 & 0 & d'\!\! \end{bmatrix*} } \\
  &\text{with} \qquad a'= \frac{-e}{b^2-ea}, \ \ b'= \frac{b}{b^2-ea}, \ \
  c' = \inv{c}, \ \ d' = \inv{d}, \ \ e'= \frac{-a}{b^2-ea}. \notag
\end{align}
In addition, $a,c,d,e>0$ (by the positive definiteness of $\brC,\brA$) and $b^2\shm ae\shneq0$ (by the invertibility of $\brC,\brA$).

Patterns~\eqref{axisym} of course apply if $\G$ is invariant under any symmetry group having $D_4$ as a sub-group. This includes the important case of rotational symmetry of $\G$ about $\brO\bfe_3$, i.e. axisymmetry. In all those cases, the pattern~\eqref{axisym} of $\brZ$ implies, on using the parametrization~\eqref{rigid:param}, that $\AUO(\bfW)=0$ for any net motion direction $\bfW$. Any $\bfW$ must in this case lie in a symmetry plane, therefore so does $\bfW^{\star}$, implying that the globally optimal motion $\bfW^{\star}$ is rotationless.

\subsection*{Three orthogonal symmetry planes (case S3)} Then, matrices $\brC,\brA,\brC^{-1},\brZ$ all are diagonal, with positive coefficients due to their positive definiteness. Being an eigenvector of $\brZUU$, the globally-optimum net motion $\bfW^{\star}$ therefore lies in the intersection of two symmetry planes, and that motion is rotationless.

\section{Axisymmetric swimmer}
\label{sec:axi}

We now address the specific case where the body shape $\G$ is axisymmetric (but not spherical), the Cartesian frame being adjusted so that $\brO\bfe_3$ is the axis of rotational symmetry. With this convention, the rigid-body basis function $\bfuR_r=\bfuR_6=\eps_{3qp}x_p\bfe_q$, which generates spinning rotations about $\brO\bfe_3$, is tangential on $\G$, so that $\bfuR_6\in\bsfHs$. In addition to exploiting the resulting matrix patterns~\eqref{axisym}, this case requires a modified version of the $r$ flow problems~\eqref{forward:PDE:fD:JL} for reasons explained next.\enlargethispage*{1ex}

\subsection{Modified solution procedure}

The intersection $\bsfR\cap\bsfHs=\text{span}(\bfuR_r)$ is now non-trivial, with the following consequences:
\begin{compactenum}[(a)]
\item The mapping $\bfuS\mapsto\bfupal[\bfuS]$, given in~\eqref{R:expr}, provides $\sfR\bfuR_r=-\bfuR_r$, implying
\begin{equation}
  \bfwS+t\bfuR_r + \sfR[\bfwS+t\bfuR_r] = \bfwS + \sfR[\bfwS]. \label{slip:ambig}
\end{equation}
In other words, applied slip velocities $\bfwS+t\bfuR_r$ produce for any $t$ the same fluid flow and swimmer motion while the added rigid-body component $\sfR[\bfwS+t\bfuR_r]$ depends on $t$ (leading to ambiguity in the definition of the net swimming motion).
\item The tractions $\bff^{\alpha}_i$, whose tangential part is the data of problem~\eqref{forward:PDE:fD:JL}, verify
\begin{equation}
  \lbra \bff^{\alpha}_i, \bfuR_r \rbra_{\G} = -\delta_{ir}. \label{nnf:violation}
\end{equation}
If $i\sheq r$, the no-net-force condition $\lbra \bff[\bfv_i], \bfuR_r \rbra_{\G} = 0$ thus cannot be enforced.
\item For each $i=1,\ldots,r$, condition $(\bfv_i\shm\bfu\Rsup_i)\tsup\bfn=0$ in problem~\eqref{forward:PDE:fD:JL} is insensitive to the $\bfuR_r$ component of $\bfu\Rsup_i$.
\end{compactenum}
As a result, the power-loss optimization method of Section~\ref{JW:opt:generic} must be modified. First, the ambiguity~\eqref{slip:ambig} is removed by seeking slip velocities in the subspace $\bsfHsr$ (of codimension 1) of $\bsfHs$ limited to spin-free tangent vector fields:
\begin{equation}
  \bsfHsr = \lcb \bfwS\in\bsfHs,\,\lbra \bfuR_r,\bfwS \rbra_{\G}=0 \rcb, \label{HTr:def}
\end{equation}
Second, we modify problems~\eqref{forward:PDE:fD:JL}, which become
\begin{equation}\left.
\begin{aligned}
  -\bfna p_i + \mu\Delta\bfv_i = \bfze, \qquad \Div\bfv_i = 0 & \qquad \text{in $\OO$} \\
  \bfPi\bff_i = \bfg\Dsup_i, \qquad (\bfv_i\shm\bfv_i\Rsup)\tsup\bfn=0 & \qquad \text{on $\G$} \\
  \lbra \bff_i,\bfuR_k \rbra_{\G} = 0 & \qquad k=1,\ldots,r\shm1 \\
  \lbra \bfv_i\shm\bfv_i\Rsup,\bfuR_r \rbra_{\G} = 0
\end{aligned}\ \right\}\qquad\qquad (1\shleq i\shleq r), \label{forward:PDE:fD:JL:axi}
\end{equation}
with the tangential traction field $\bfg\Dsup_i$ given by
\begin{equation}
  \bfg\Dsup_i = \bfPi\bff^{\alpha}_i \ \ (i\shneq r), \qquad \bfg\Dsup_r = \bfPi\bff^{\alpha}_r + \inv{\lbra \bfuR_r,\bfuR_r \rbra_\Gamma}\bfuR_r , \label{fD:galerkin:JL:axi}
\end{equation}
with $\bff^{\alpha}_i$ as in~\eqref{R:expr}; this definition is designed to satisfy $\lbra \bfg\Dsup_r,\bfuR_r \rbra_{\G}=0$ (see~\eqref{nnf:violation}), i.e. to fulfill the $r$-th no-net-force condition. The solutions of problems~\eqref{forward:PDE:fD:JL:axi} allows to define the slip velocities
\begin{equation}
  \bfz_i := (\bfv_i\shm\bfv_i\Rsup)|_{\G} \in \bsfHsr, \label{zi:def:JL:axi}
\end{equation}
which, by verifying
\begin{equation}
  \lbra \bff[\bfv_i],\bfwS \rbra_{\G} = \lbra \bff^{\alpha}_i,\bfwS \rbra_{\G} = \alpha_i[\bfwS],
\end{equation}
extract the rigid-body components induced by any slip velocity $\bfwS\shin\bsfHsr$ as per~\eqref{R:expr}. Then, definitions~\eqref{eig:rank:r:JL} of the matrix $\brA$ and~\eqref{y:def} of the slip velocities $\bry$ apply as in the general case. We thus set as before $\bsfH_r=\text{span}(\bfy_1,\ldots,\bfy_r)$; $\bsfH_r$ is a $r$-dimensional subspace of $\bsfHsr$ which produces spin-free slip velocities achieving least power loss for given $\bfU,\bfOm$.

\subsection{Power loss minimization}
\label{sec:algo:axi}

Again invoking the parametrization~\eqref{rigid:param} of $\bfU,\bfOm$, the pattern~\eqref{axisym} of $\brZ$ implies $\AUO(\bfW)=0$ for any net direction $\bfW$, so that~\eqref{s=0} applies; we obtain
\begin{equation}
\begin{gathered}
  \bfOm=\bfze, \qquad
  \AUU = \inv{Z_{11}}(W_1^2\shp W_2^2) + \inv{Z_{33}}W_3^2, \\
  \bfU = \inv{\AUU}\Lpar \frac{W_1}{Z_{11}},\,\frac{W_2}{Z_{11}},\,\frac{W_3}{Z_{33}} \Rpar, \qquad \PcalH(\bfW)=\inv{\AUU}.
\end{gathered}\label{s=0:axi}
\end{equation}
In the exceptional cases where $\G$ is such that $Z_{11}\sheq Z_{33}$, we have $\bfU\sheq \bfW$ and $\PcalH(\bfW)\sheq Z_{11}$ for any direction $\bfW$: the same (optimal) power loss is incurred by any slip velocity~\eqref{y:def} with $\bfU\sheq \bfW,\bfOm\sheq \bfze$. When $Z_{11}\shneq Z_{33}$, $\bfW$ must be either aligned with, or orthogonal to, the rotation axis to produce valid net motions ($\bfU\sheq \bfW$), entailing respective power losses $\PcalH(\bfW)\sheq Z_{33}$ or $\PcalH(\bfW)\sheq Z_{11}$; the direction producing $\PcalH(\bfW)\sheq \text{min}(Z_{11},Z_{33})$ then defines the globally optimal motion. Any such optimal motion is straight (i.e. rotationless), corroborating the analysis of~\cite{stone:96}. The computational treatment for a given swimmer thus relies on the evaluation of only $Z_{11},Z_{33},Z_{15}$ and the tangential fields $\bfy_1,\bfy_3$, and requires solving 6 (instead of 12) flow problems overall. This results in the following procedure, where steps (i)-(iv) adapt the corresponding steps of Sec.~\ref{sec:algo:generic}, and which requires only six flow solutions (instead of 12 for general 3D shapes) and four selected entries of both $\brC$ and $\brA$.\medskip

\begin{compactenum}[(i)]
\item Solve problems~\eqref{Stokes}-\eqref{BC} with $\bfuD=\bfuR_{\ell}$ for $\ell=1,3,4$, obtain corresponding tractions $\bff\Rsup_{\ell}$.
\item Compute entries $C_{11}\sheq\lbra \bffR_1,\bfuR_1 \rbra_{\G}$, $C_{33}\sheq\lbra \bffR_3,\bfuR_3 \rbra_{\G}$,  $C_{44}\sheq\lbra \bffR_4,\bfuR_4 \rbra_{\G}$, $C_{15}\sheq\lbra \bffR_1,\brr\bfuR_4(\brr^{-1}\dotp) \rbra_{\G}$ of $\brC$.\\
Evaluate $D_{\rmC}\sheq C_{15}^2\shm C_{11}C_{44}$. Set up rigid body extraction tractions $\bff^{\alpha}_1\sheq(C_{11}\bffR_1-C_{15}\brr\bffR_4)/D_{\rmC}$, $\bff^{\alpha}_3\sheq-\bffR_3/C_{33}$, $\bff^{\alpha}_4\sheq(C_{44}\bffR_4+C_{15}\brr\bffR_1)/D_{\rmC}$.
\item Solve problems~\eqref{forward:PDE:fD:JL:axi} for $i\sheq1,3,4$, obtain corresponding $\bfz_i$ given by~\eqref{zi:def:JL:axi} and tractions $\bff_i$.
\item Compute entries $A_{11}\sheq\lbra \bff_1,\bfz_1 \rbra_{\G}$, $A_{33}\sheq\lbra \bff_3,\bfz_3 \rbra_{\G}$,  $A_{44}\sheq\lbra \bff_4,\bfz_4 \rbra_{\G}$, $A_{15}\sheq\lbra \bff_1,\brr\bfz_4(\brr^{-1}\dotp) \rbra_{\G}$ of $\brA$. Evaluate $D_{\rmA}\sheq A_{15}^2\shm A_{11}A_{44}$; compute $Z_{11}\sheq-A_{44}/D_{\rmA}$, $Z_{33}\sheq1/A_{33}$ and $Z_{15}\sheq A_{15}/D_{\rmA}$.\\
Compute the tangential fields $\bfy_1\sheq Z_{11}\bfz_1+Z_{15}\brr\bfz_4$ and $\bfy_3\sheq Z_{33}\bfz_3$.
\item If $Z_{33}\leq Z_{11}$, set $\bfW^{\star}=\bfe_3$, $\Pcal^{\star}=Z_{33}$ and $\bfuS{}^{\star}=\bfy_3$.\\
If $Z_{11}< Z_{33}$, set $\bfW^{\star}=\bfe_1$, $\Pcal^{\star}=Z_{11}$ and $\bfuS{}^{\star}=\bfy_1$.\medskip
\end{compactenum}
Problems~\eqref{forward:PDE:fD:JL:axi} all needed to be formulated in order to define the matrix $\brA$ and tangential fields $\brz$, $\bry$, so that the pattern~\eqref{axisym} of $\brZ$ and the rotationless character of optimal motions coud be inferred. In step (v), if the swimmer verifies $Z_{11}=Z_{33}$, the choice $\bfW=\bfe_3$ (net motion along the axis of rotational symmetry) is arbitrary, as the same optimal power loss can be achieved for any choice of $\bfW$.

\begin{figure}[t] \centering
  \includegraphics[trim=100 0 190 250, clip, width=0.78\textwidth]{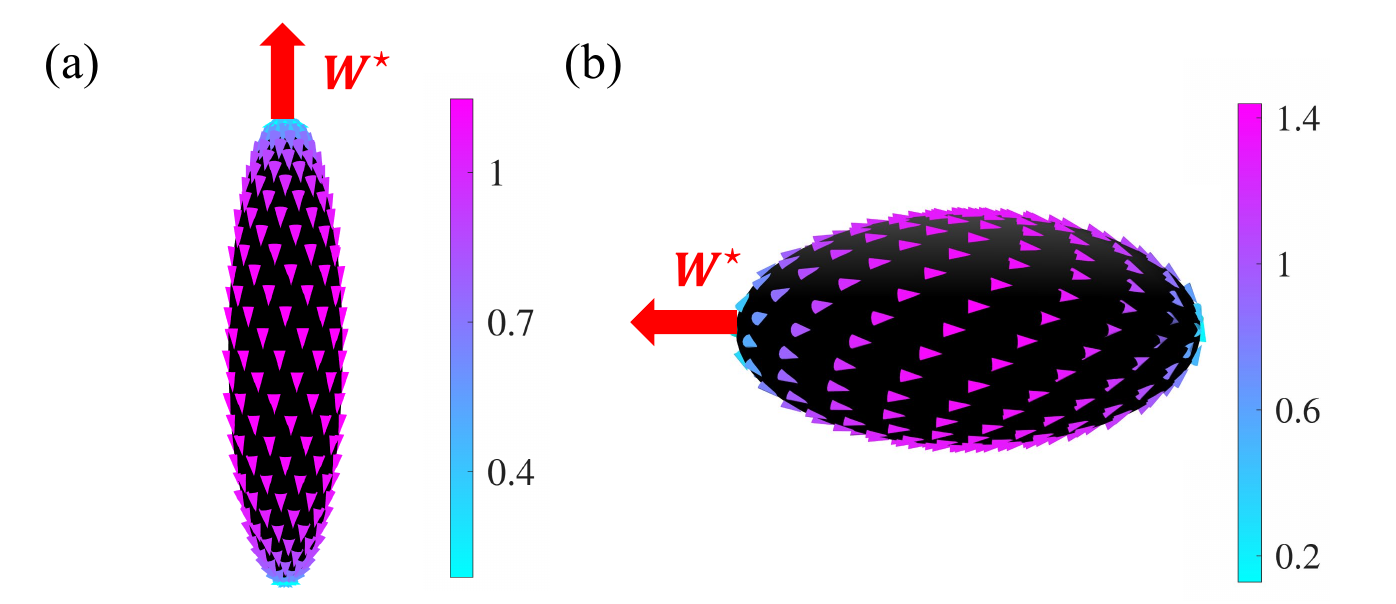}
  \caption{Two axisymmetric shapes (one swimmer very elongated along the axial symmetry direction, the other almost flat) producing globally optimal motions that are respectively axial and transverse (i.e. in the plane orthogonal to the axis). (a) Prolate spheroid with semimajor axes 1/4 ,1/4, 1. (b) Oblate spheroid with semimajor axes 1, 1,1/2.}\label{fig:shapeiterations3}
\end{figure}

Figure~\ref{fig:shapeiterations3} illustrates the finding that globally optimal motions of axisymmetric swimmers are either axial $(Z_{33}\shl Z_{11})$ or transverse $(Z_{11}\shl Z_{33})$ depending on the axial elongation of the swimmer.\enlargethispage*{3ex}

\section{Conclusions}
\label{sec:conclu}

We have presented a rigorous mathematical framework for the slip optimization of 
microswimmers with arbitrary three-dimensional geometries in Stokes flow. The central 
contribution is a dimension reduction result (Proposition~\ref{y:props} and 
Lemma~\ref{lemma:z}) showing that the infinite-dimensional power loss minimization 
problem~\eqref{JL:opt:generic} is exactly equivalent to a low-dimensional programming problem posed over a specific $r$-dimensional subspace $\bsfH_r$ of tangential 
slip velocities, where $r = 6$ for general three-dimensional shapes. The basis functions 
spanning $\bsfH_r$ are defined from the solutions of $2r$ auxiliary Stokes flow problems, 
and have the property that any $\bfuS \in \bsfH_r$ achieves strictly less power loss than 
any slip velocity outside $\bsfH_r$ producing the same rigid-body motion. This allows 
the optimal slip profile to be determined without repeated flow solves in the optimization 
loop; once the system matrices $\brC$ and $\brA$ are assembled, the remaining optimization 
is purely algebraic.\enlargethispage*{1ex}

The optimality conditions for the power loss minimization were analyzed in detail 
(Section~\ref{sec:solv:opt}). When the direction of net motion $\bfW$ is prescribed, 
the partial minimization problem~\eqref{JL:opt:rdim:part} has a closed-form solution 
(Lemma~\ref{partial:min}), giving the optimal spin parameter $s$, drift velocity $\bfV$, 
and rigid-body parameters explicitly in terms of entries of the matrix $\brZ = \brA^{-1}$. 
The reduced global minimization over $\bfW \in \Sigma$ does not generally admit a closed 
form, but any of its solutions satisfying $\AUO(\bfW) = 0$ must be an eigenvector of 
$\brZUU^{-1}$ (Lemma~\ref{AUO:opt}), ensuring that globally optimal motions are always 
physically consistent in the sense of~\eqref{W:def}.

A systematic analysis of the influence of swimmer symmetry on the structure of the 
matrices $\brC$ and $\brA$ (Section~\ref{CA:sym}) shows that reflectional and rotational 
symmetries induce block-sparsity patterns in these matrices and their inverses. 
Proposition~\ref{prop:optim:symm} identifies two regimes: for swimmers with axisymmetry 
or three orthogonal symmetry planes, the globally optimal motion must be 
rotationless; for swimmers with fewer symmetries, rotational globally optimal motions are 
possible, and are indeed observed numerically. These structural results reduce the 
computational cost in the axisymmetric case to $6$ flow solutions and the evaluation 
of only four entries of $\brC$ and $\brA$, rather than the $2r = 12$ solutions required 
for a general three-dimensional shape.

The numerical treatment relies on a high-order boundary integral method to discretize 
the auxiliary flow problems~\eqref{Stokes}-\eqref{BC} and~\eqref{forward:PDE:fD:JL}. Spectral 
accuracy is demonstrated in Figure~\ref{fig: convergence}, and the optimal slip profiles 
for axisymmetric shapes are shown in our companion paper~\cite{B-2025-05} to recover the results of~\cite{B-2020-04} 
to high precision. Numerical results for non-convex and chiral shapes confirm the 
theoretical predictions regarding the role of symmetry in determining whether globally 
optimal motions are rotational or translational.

Several directions remain open. First, the present 
framework addresses steady slip velocities and the minimization of hydrodynamic power 
loss; extending the theory to time-periodic slip profiles, as required to capture 
metachronal wave patterns, or to the joint optimization of swimmer shape and slip, as 
considered in the axisymmetric setting in~\cite{B-2023-10}, would broaden the range of 
accessible problems considerably. Second, the connection between the matrix $\brA^{-1} = 
\brZ$ and the classical resistance and mobility tensors of low-Reynolds-number 
hydrodynamics deserves further investigation; in particular, clarifying how the 
partitions $\brZUU$, $\brZUO$, $\brZOO$ relate to known tensor structures for bodies 
of special geometry could illuminate the physical interpretation of the optimal 
spin parameter $s$ and drift velocity $\bfV$. Finally, incorporating hydrodynamic 
interactions in suspensions of multiple swimmers, or the influence of confining boundaries, 
within the present variational framework would be a natural and practically significant 
extension.

\section*{Acknowledgements}
KD acknowledges support from NSF GRFP under grant DGE-2241144. SV acknowledges support from NSF under grant DMS-2513346.
\appendix

\section{Proof of Lemma~\ref{recipr}}
\label{proof:recipr}
Let $(\bfu_i,p_i,\bff_i)$ solve problem~\eqref{Stokes}-\eqref{BC} with Dirichlet datum $\bfuD_i$ ($i=1,2$). Taking the inner product of the first of~\eqref{Stokes} obeyed by $(\bfu_1,p_1)$ with $\bfu_2$ and invoking Green's identity and incompressibility, we find
\begin{equation}
  \iO 2\mu\bfD[\bfu_1]\dip\bfD[\bfu_2] \dV = \lbra\bff_1,\bfu_2 \rbra_{\G}. \label{forward:weak}
\end{equation}
The variational identity~\eqref{forward:weak} is valid for the exterior problem~\eqref{Stokes}-\eqref{BC} since $\bfu(\bfx)\sheq O(|\bfx|^{-1})$ and $p(\bfx)\sheq O(|\bfx|^{-2})$ as $|\bfx|\to\infty$ (e.g. by virtue of their respective integral representations), ensuring the boundedness of both terms in~\eqref{forward:weak}. The first identity of Lemma~\ref{recipr} then corresponds to $\bfuD_1=\bfuD_2=\bfuD$, while the second results from the symmetry in $\bfu_1,\bfu_2$ of the left-hand side of~\eqref{forward:weak}.

\section{Proof of Lemma~\ref{lemma:z}}
\label{lemma:z:proof}
For item (i), the complete velocity on $\G$ is given by~\eqref{zi:def:JL} as $\bfv_i=\bfz_i\shp\bfv_i\Rsup$. By classical uniqueness results for Stokes flows, the full solution of problem~\eqref{forward:PDE:fD:JL} thus also solves problem~\eqref{Stokes}-\eqref{BC} with $\bfuD=\bfz_i\shp\bfv_i\Rsup$. That $\bfv\Rsup_i=\sfR\bfz_i$ follows directly from the no-net-force conditions~\eqref{no-net-force} being prescribed in problems~\eqref{forward:PDE:fD:JL} and the definitions of the other quantities.

Regarding item (ii), we have
\begin{equation}
  \lbra \bff[\bfz_i\shp\bfv\Rsup_i],\bfw \rbra_{\G}
  = \lbra  \bfPi\bff[\bfv_i],\bfw \rbra_{\G}
  = \lbra \bff^{\alpha}_i,\bfw \rbra_{\G} = \alpha_i[\bfw]
\end{equation}
for any $\bfw\in\bsfHs$, by virtue of~\eqref{R:expr}.

Turning finally to item (iii), any $\bfuS_1\shin\bsfH_r$ has the form $\bfuS_1=\brz\bfal$ for some $\bfal\shin\Rbb^r$, so that we have
\begin{equation}
  2\lbra \bff[\bfuS_1\shp\sfR\bfuS_1],\bfuS_2 \rbra_{\G}
  = 2\bfal\tsup\bff[\brz\tsup\shp\sfR\brz\tsup],\bfuS_2 \rbra_{\G} = 2\bfal\tsup\bfupal[\bfuS_2] = 0
\end{equation}
which, used in~\eqref{Jw:expand}, gives the clained inequality. The proof of the lemma is complete.

\section{Proof of Proposition~\ref{y:props}}
\label{sec:proof:y:props}
For item (i), let $\bfal\shin\Rbb^r$. By property~\eqref{FR:yS}, the slip velocity $\bfuS=\bry\bfal\in\bsfH_r$ verifies $\bfupal[\bfuS]=\bfal$, so that the $\bsfH_r\to\Rbb^r$ mapping $\bfuS\mapsto\bfupal[\bfuS]$ is surjective. Since $\bsfH_r$ and $\Rbb^r$ are finite-dimensional with the same dimension, the mapping is bijective by the rank-nullity theorem.

Item (ii) follows from evaluating $P(\bfuS)\sheq\lbra \bff[\bfuS\shp\sfR\bfuS],\bfuS \rbra_{\G}$ by writing $\bfuS\sheq\brz\brA^{-1}\bfal$ (using the definition~\eqref{y:def} of $\bry$ and~\eqref{eig:rank:r:JL}). Finally, item (iii) just restates Lemma~\ref{lemma:z}(iii).

\section{Proof of Lemma~\ref{partial:min}}
\label{sec:proof:partial:min}
To accommodate the equality constraint in problem~\eqref{JL:opt:rdim:part}, we define the Lagrangian
\begin{equation}
  L(\bfV,s,\xi;\bfW)
 := \Pcal(s,\bfV,\bfW) - 2\xi\bfV\tsup\bfW, \label{aux24}
\end{equation}
whose first-order stationarity conditions (a) $\del[\xi]L=0$, (b) $\del[\bfV]L=0$ and (c) $\del[s]L=0$ for $L$ read
\begin{equation}
  \text{(a) \ } \bfV\tsup\bfW = 0, \qquad
\begin{aligned}
  \text{(b) \ }\brZUU\bfV
  + \lpar s\brZUO + \brZUU \rpar\bfW - \xi\bfW &= \bfze, \\
  \text{(c) \ } \bfV\tsup \brZUO \bfW
    + \bfW\tsup\lpar s\brZOO + \brZUO \rpar \bfW &= 0.
\end{aligned} \label{aux22}
\end{equation}
The combinations $\bfW\tsup\brZOU\brZUU^{-1}$(b)-(c) and $\bfW\tsup\brZUU^{-1}$(b)-(a)  eliminate $\bfV$ and give
\begin{equation}
  \AOO s + \AUO\xi = 0, \qquad -\AUO s + \AUU\xi = |\bfW|^2,
\end{equation}
with $\AUU,\,\AUO,\,\AOO$ as defined in~\eqref{aux21}. Solving the above produces $s,\xi$ as given by~\eqref{s:mu:expr}. Then, using (b) in $\Pcal(s,\bfV,\bfW)$ given by~\eqref{aux20} yields, with the help of the alternative expression
\begin{align}
  \Pcal(s,\bfV,\bfW)
 &= \lsqb \brZUU\bfV + \lpar s\brZUO + \brZUU \rpar\bfW \rsqb\tsup \brZUU^{-1}
    \lsqb \brZUU\bfV + \lpar s\brZUO + \brZUU \rpar\bfW \rsqb \suite\qquad
    + s^2\bfW\tsup \lsqb \brZOO - \brZOU\brZUU^{-1}\brZUO \rsqb \bfW,
\end{align}
the minimum value of the power loss as
\begin{equation}
  \PcalH(\bfW)
 = \AUU\xi^2 + \AOO s^2
 = \frac{|\bfW|^4\AOO}{D} = \frac{|\bfW|^4}{\AUU} - \frac{|\bfW|^4(\AUO)^2}{\AUU D} = |\bfW|^2\xi.
\end{equation}
Finally, the optimal values of $\bfOm$ and $\bfU$ given $\bfW$ are obtained from~\eqref{rigid:param}, with $\bfV$ found from~(\ref{aux22}b) and $s$ as given by~\eqref{s:mu:expr}.

\section{Proof of Lemma~\ref{AUO:opt}}
\label{AUO:opt:proof}
Any $\bfW\in\Sigma$ solving problem~\eqref{JW:min:angular} must satisfy the first-order conditions
\begin{equation}
   \PcalH'(\bfW;\bfz) - 2\nu\bfz\tsup\bfW = 0 , \qquad |\bfW|^2=1 \qquad \text{for all }\bfz\shin\Rbb^3,
\label{JW:1storder}
\end{equation}
where $\nu$ is the Lagrange multiplier associated with the unit-norm equality constraint $\bfW\in\Sigma$ and $\PcalH'(\bfW;\bfz)$, the directional derivative of $\PcalH$ at $\bfW$ in the direction $\bfz$, is found to be given by
\begin{align}
  \MoveEqLeft[8]{
  \PcalH'(\bfW;\bfz)
 = \frac{|\bfW|^4}{D^2(\bfW)} \Lsqb A^2_{U\OO}(\bfW) A'_{\OO\OO}(\bfW;\bfz)
 - A^2_{\OO\OO}(\bfW) A'_{UU}(\bfW;\bfz) } \notag \\[-1ex]
 & - 2\AOO(\bfW)\AUO(\bfW) A'_{U\OO}(\bfW;\bfz) \Rsqb
 + \frac{4|\bfW|^2\AOO(\bfW)}{D(\bfW)}\,(\bfz\tsup\bfW)
\end{align}
If such $\bfW$ also verifies $\AUO(\bfW)=0$, the first stationarity condition in~\eqref{JW:1storder} becomes
\begin{equation}
  \Lpar \frac{4|\bfW|^2}{\AUU(\bfW)} - 2\nu \Rpar\,(\bfz\tsup\bfW)
 -\frac{|\bfW|^4\,A'_{UU}(\bfW;\bfz)}{A^2_{UU}(\bfW)} = 0 \qquad \text{for all }\bfz\shin\Rbb^3
\end{equation}
which, by virtue of $A'_{UU}(\bfW;\bfz)=2\bfz\tsup\bfZ^{-1}_{UU}\bfW$, yields
\begin{equation}
  |\bfW|^4\bfZ^{-1}_{UU}\bfW + \lsqb A^2_{UU}(\bfW)\nu - 2|\bfW|^2A_{UU}(\bfW) \rsqb \bfW = 0.
\end{equation}
To satisfy the above (homogeneous in $\bfW$) condition, $\bfW$ must be an eigenvector of $\bfZ^{-1}_{UU}$, as claimed.\enlargethispage*{1ex}

\bibliographystyle{siamplain}%
\bibliography{marcbibs, refs, refs2}%
\end{document}